\documentclass{amsart}
\usepackage{enumerate, vmargin}
\usepackage[arrow,matrix]{xy}

\newtheorem{thm}{Theorem}[section]
\newtheorem{prop}[thm]{Proposition}
\newtheorem{cor}[thm]{Corollary}
\newtheorem{lem}[thm]{Lemma}
\newtheorem{defi}[thm]{Definition}
\newtheorem{remark}[thm]{Remark}
\newtheorem{example}[thm]{Example}
\newtheorem{pb}[thm]{Problem}

\newenvironment{rk}{\begin{remark}\rm}{\end{remark}}

\newcommand{\real}{{\mathbb R}}
\newcommand{\nat}{{\mathbb N}}
\newcommand{\ent}{{\mathbb Z}}
\newcommand{\comp}{{\mathbb C}}

\newcommand{\F}{{\mathcal F}}
\newcommand{\G}{{\Gamma}}
\renewcommand{\H}{{\mathcal H}}
\newcommand{\M}{{\mathcal M}}
\newcommand{\U}{{\mathcal U}}

\renewcommand{\a}{\alpha}
\renewcommand{\b}{\beta}
\newcommand{\g}{\gamma}
\renewcommand{\d}{\delta}

\newcommand{\e}{\varepsilon}
\renewcommand{\th}{\theta}
\renewcommand{\l}{\lambda}
\newcommand{\La}{\Lambda}
\newcommand{\f}{\varphi}
\renewcommand{\O}{{\Omega}}
\renewcommand{\o}{{\omega}}
\newcommand{\s}{\sigma}

\newcommand{\el}{\ell}

\newcommand{\Tr}{\mbox{\rm Tr}}
\newcommand{\tr}{\mbox{\rm tr}}
\newcommand{\ot}{\otimes}
\newcommand{\op}{\oplus}

\newcommand{\la}{\langle}
\newcommand{\ra}{\rangle}
\newcommand{\wt}{\widetilde}
\newcommand{\wh}{\widehat}
\newcommand{\n}{\noindent}
\newcommand{\8}{\infty}

\newcommand{\pf}{\noindent{\it Proof.~~}}
\newcommand{\cqd}{\hfill$\Box$}
\newcommand{\be}{\begin{eqnarray*}}
\newcommand{\ee}{\end{eqnarray*}}
\newcommand{\beq}{\begin{equation}}
\newcommand{\eeq}{\end{equation}}

\begin{document}

\title[Operator space Grothendieck inequalities]
{Operator space Grothendieck inequalities for
noncommutative $L_p$-spaces}
\author{Quanhua XU}
\address{Laboratoire de Math\'{e}matiques\\
Universit\'{e} de France-Comt\'{e}\\
16 Route de Gray\\
25030 Besan\c con Cedex,  France } \email{qx@math.univ-fcomte.fr}

\numberwithin{equation}{section}


\begin{abstract}
 We prove the operator space Grothendieck inequality for bilinear
forms on subspaces of  noncommutative $L_p$-spaces with $2<p<\8$.
One of our results states that given a map  $u: E\to F^*$, where
$E, F\subset L_p(M)$ ($2<p<\8$, $M$ being a von Neumann algebra),
$u$ is completely bounded iff $u$ factors through a direct sum of
a $p$-column space and a $p$-row space. We also obtain several
operator space versions of the classical little Grothendieck
inequality for maps defined on a subspace of a noncommutative
$L_p$-space ($2<p<\8$) with values in a $q$-column space for every
$q\in [p', p]$ ($p'$ being the index conjugate to $p$). These
results are the $L_p$-space analogues of the recent works on the
operator space Grothendieck theorems by Pisier and Shlyakhtenko.
The key ingredient of our arguments is some Khintchine type
inequalities for Shlyakhtenko's generalized circular systems. One
of our main tools is a Haagerup type tensor norm, which turns out
particularly fruitful when applied to subspaces of noncommutative
$L_p$-spaces ($2<p<\8$). In particular, we show that the norm dual
to this tensor norm, when restricted to subspaces of
noncommutative $L_p$-spaces, is equal to the factorization norm
through a $p$-row space.
\end{abstract}

 \setcounter{page}{0}
 \setcounter{section}{-1}
 \maketitle
 \thispagestyle{empty}


\vskip 3cm

\n{\bf Plan:}
\begin{itemize}
\item[0.]  Introduction
\item[1.]  Preliminaries
\item[2.]  A Haagerup type tensor norm and factorization through
 $p$-row spaces
\item[3.] Noncommutative Khintchine inequalities
\item[4.] Vector-valued noncommutative $L_p$
\item[4.]  Proof of Theorem \ref{G}
\item[6.]  Proof of Theorem \ref{LG}
\item[7.] Applications
\end{itemize}

\newpage


\section{ Introduction}\label{intro}

In the remarkable recent work \cite{pisshlyak}, Pisier and
Shlyakhtenko obtained the operator space version of the famous
Grothendieck theorem. This can be stated as follows. Let $E,
F\subset B(H)$ be operator spaces and $u: E\times F\to\comp$ a
bilinear form. Assume $E$ and $F$ exact. Then $u$ is jointly
completely bounded iff there are a constant $K$ and states $f_i,
g_i$ ($i=1,2$) on $B(H)$ such that
 $$
 |u(a, b)|\le K\Big[\big(f_1(aa^*)\,g_1(b^*b)\big)^{1/2}+
 \big(f_2(a^*a)\,g_2(bb^*)\big)^{1/2}\Big],\quad a\in E,\; b\in F.
 $$
Moreover, if $K$ denotes the least constant in the inequality
above, then $K\approx \|u\|_{jcb}$, where the relevant equivalence
constants depend only on the exactness constants of $E$ and $F$.
We refer to the next section for background on operator space
theory and all unexplained notions. \cite{pisshlyak} also contains
several interesting variants of the above statement, especially
when both $E$ and $F$ are C*-algebras (then the exactness
assumption is needed for only one of them). In this latter case,
the corresponding inequality is exactly the version for operator
space theory of the noncommutative Grothendieck inequality
obtained first by Pisier \cite{pis-groC*} with an approximability
assumption, and then in the full generality by Haagerup
\cite{haag-gro}. Recall that the classical Grothendieck inequality
corresponds to the case where both $E$ and $F$ are commutative
C*-algebras (in the Banach space theory). \cite{pis-cbms} is an
excellent reference for the classical and noncommutative
Grothendieck inequalities for Banach spaces.

On the other hand, Maurey \cite{mau-asterisque} extended the
classical Grothendieck inequality to bilinear forms on commutative
$L_p$-spaces. (In this paper, we use ``commutative $L_p$-spaces''
to distinguish the usual $L_p$-spaces from the general
noncommutative $L_p$-spaces, which are our main objects.) Let
$(\O, \mu)$ be a measure space and $2<p, q\le\8$. Let $u: L_p(\O,
\mu)\times L_q(\O, \mu)\to\comp$ be a bilinear form. Then $u$ is
bounded iff there are a constant $K$ and positive unit functionals
$f\in \big(L_{p/2}(\O,\mu)\big)^*$, $g\in
\big(L_{q/2}(\O,\mu)\big)^*$ such that
 $$
 |u(a, b)|\le K\big[f(|a|^2)\,g(|b|^2)\big]^{1/2},
 \quad a\in  L_p(\O, \mu),\; b\in  L_q(\O, \mu).
 $$
Again the best constant $K$ is equivalent to $\|u\|$. Although
this statement is not explicitly stated in \cite{mau-asterisque},
it immediately follows from Kwapien's theorem (cf. \cite
[Corollary 3.6] {pis-cbms}) and the little Grothendieck theorem
for $L_p$-spaces in [M1]. This latter theorem says that if $u:
L_p(\O, \mu)\to H$ is a bounded map ($2<p\le\8$; $H$ being a
Hilbert space), then there is a positive unit functional  $f\in
\big(L_{p/2}(\O,\mu)\big)^*$ such that
 $$|u(a)|\le K_0\|u\|\, \big[f (|a|^2)\big]^{1/2},\quad
 a\in L_p(\O, \mu),$$
where $K_0$ is a universal constant.

It is this last statement which was extended to the noncommutative
setting by Lust-Piquard \cite{lust-gro}. Let $M$ be a semifinite
von Neumann algebra equipped with a normal semifinite faithful
trace $\tau$. Let $L_p(M)$ be the associated noncommutative
$L_p$-space. Let $2<p\le\8$ and $u: L_p(M)\to H$ be a bounded map.
Then Lust-Piqard's theorem claims that there are positive unit
functionals $f_1, f_2\in \big(L_{p/2}(M)\big)^*$ such that
 \beq\label{piquard}
 |u(a)|\le K_0\|u\|\,\big[f_1(aa^*)+f_2(a^*a)\big]^{1/2},
 \quad a\in  L_p(M).
 \eeq
As in the commutative case, using Kwapien's theorem, we then
deduce that for any bounded bilinear form $u: L_p(M)\times L_q(M)
\to\comp$ ($2<p, q\le\8$) there are positive unit functionals
$f_i\in \big(L_{p/2}(M)\big)^*$ and  $g_i\in
\big(L_{q/2}(M)\big)^*$ such that
 \beq\label{piquard1}
 |u(a, b)|\le K\|u\|\, \big[f_1(aa^*)+f_2(a^*a)\big]^{1/2}
 \,\big[g_1(b^*b)+g_2(bb^*)\big]^{1/2} ,\
 a,b\in L_p(M).
 \eeq

\medskip

Having all these in mind, one is naturally attempted to find out
the operator space versions  of (0.1) and (0.2) in the spirit of
\cite{pisshlyak} for bilinear forms on noncommutative
$L_p$-spaces. This is the main concern of the present paper. The
noncommutative $L_p$-spaces we use are those constructed by
Haagerup \cite{haag-Lp}. Thus type III von Neumann algebras are
also allowed. The reader is referred to the next section for a
brief introduction to noncommutative $L_p$-spaces and their
operator space structure. To state our main result, we need  some
notations (see the next section for more details). Given a Hilbert
space $H$ and $1\le p\le\8$ we denote by $H_p^c$ (resp. $H_p^r$)
the Schatten $p$-class $S_p(\comp, H)$ (resp. $S_p(\bar H,
\comp)$) equipped with its natural operator space structure.
$H_p^c$ (resp. $H_p^r$) is the column (resp. row) subspace of the
Schatten class $S_p(H)$ (resp. $S_p(\bar H)$). When $H$ is
separable and infinite dimensional, $H_p^c$ and $H_p^r$ are
respectively denoted by $C_p$ and $R_p$. Given operator spaces $X,
E$ and $F$ we denote by $\G_X(E, F)$ the family of all maps $u:
E\to F$ which factors through $X$ by c.b. maps. We will need this
notably when $X$ is one of $R_p, C_p$ and $R_p\op_p C_p$. If $u:
E\times F\to\comp$ is a bounded bilinear form, $\wt u: E\to F^*$
denotes the associated linear map, i.e. $\wt u(x)(y)=u(x,y)$ for
all $x\in E$ and $y\in F$.

\medskip

The following is the  Grothendieck theorem for noncommutative
$L_p$-spaces in the category of operator spaces.

\begin{thm}\label{G}
 Let $M$ be a von Neumann algebra and $2< p<\8$. Let
$E, F\subset L_p(M)$ be two closed subspaces. Let $u: E\times F\to
\comp$ be a bilinear form. Then the following assertions are
equivalent
 \begin{enumerate}[{\rm i)}]
 \item  $u$ is jointly completely bounded and $\|u\|_{jcb}\le K_1$.
 \item  For any finite sequences $(a_k)\subset E$, $(b_k)\subset
F$ and $(\mu_k)\subset \real_+$
 \begin{eqnarray}\label{G2}
 \begin{array}{ccl}
 \begin{dis}\big|\sum_ku(a_k, b_k)\big|
 \end{dis}
 &\le& \begin{dis} K_2
 \Big[\big\|\big(\sum_k\mu_k a_k^*a_k\big)^{1/2}\big\|_p
 + \big\|\big(\sum_k\mu_k^{-1}a_ka_k^*\big)^{1/2}\big\|_p\Big]
 \end{dis}\\
 &&\begin{dis}\hskip 0.2cm \bullet\,\Big[\big\|
 \big(\sum_k\mu_k b_k^*b_k\big)^{1/2}\big\|_p
 + \big\|\big(\sum_k\mu_k^{-1}b_kb_k^*\big)^{1/2}\big\|_p\Big]
 \end{dis}.
 \end{array}
 \end{eqnarray}
 \item  There are positive functionals $f_1, f_2$ and $g_1, g_2$
in the unit ball of $\big(L_{p/2}(M)\big)^*$  such that for any
$(a, b)\in E\times F$
 \beq\label{G3}
 |u(a, b)|\le K_3\Big[\big(f_1(aa^*)\,g_1(b^*b)\big)^{1/2} +
 \big(f_2(a^*a)\,g_2(bb^*)\big)^{1/2}\Big].
 \eeq
 \item  For any finite sequences $(a_k)\subset E$, $(b_k)\subset
F$ and $(\mu_k)\subset \real_+$
 \begin{eqnarray}\label{G4}
 \begin{array}{ccl}
 \begin{displaystyle}
 \big|\sum_ku(a_k, b_k)\big|
 \end{displaystyle}
 &\le& \begin{displaystyle}K_4
 \Big[\big\|\big(\sum_k a_k a_k^*\big)^{1/2}\big\|_p\;
 \big\|\big(\sum_kb_k^*b_k\big)^{1/2}\big\|_p
 \end{displaystyle}\\
 &&\hskip 0.4cm  +\begin{displaystyle}\big\|\big(\sum_k\mu_ka_k^*a_k
 \big)^{1/2}\big\|_p\;
 \big\|\big(\sum_k\mu_k^{-1}b_kb_k^*\big)^{1/2}\big\|_p\Big].
 \end{displaystyle}
 \end{array}
 \end{eqnarray}
 \item  $u$ admits a decomposition $u=u_1+ u_2$, where
$u_1$ and $u_2$ are bilinear forms on $E\times F$ such that the
associated linear maps $\widetilde u_1, \ \widetilde u_2: E\to
F^*$ satisfy
  $$\widetilde u_1\in \G_{R_p}(E, F^*),\
  \widetilde u_2\in \G_{C_p}(E, F^*)\quad \mbox{and}\quad
  \max\big\{\g_{R_p}(\widetilde u_1), \g_{C_p}(\widetilde u_2)
  \big\}\le K_5.$$
 \item  $\widetilde u\in \G_{R_p\oplus_p C_p}(E, F^*)$ and
$\g_{R_p\oplus_p C_p}(\widetilde u)\le K_6$.
 \end{enumerate}
 Here the $K_i$ are constants; moreover, the best $K_i$
are equivalent uniformly in $p$, $E$ and $F$, i.e. there is an
absolute positive constant $c$ such that $c^{-1}K_j^{-1}\le K_i\le
cK_j$ for all $i, j=1,..., 6$.
 \end{thm}

This theorem  is the $L_p$-space version of the corresponding
results in \cite{pisshlyak} (see Theorems~0.3, 0.4, 0.5 and
Corollary 0.7 there). More precisely, Pisier-Shlyakhtenko's
results correspond to Theorem~0.1 above in the case of $p=\8$. In
this case, one needs the exactness assumption on $E$ and $F$,
namely, either both $E$ and $F$ are exact, or both $E$ and $F$ are
C*-algebras and one of them is exact.

\medskip

From Theorem~0.1 we can easily deduce the version of (0.1) for
operator spaces, namely, the noncommutative little Grothendieck
theorem in the category of operator spaces (see section
\ref{appli}). However, we will prove a more general result in the
spirit of \cite{pis-Rp}. This is the following theorem. Recall
that the main result of \cite{pis-Rp} corresponds again to the
case $p=\8$. Then $E$ must be supposed to be either exact or a
C*-algebra.

\begin{thm}\label{LG}
 Let $E\subset L_p(M)$ be a subspace with
$2<p<\8$ and $H$ a Hilbert space. Let $0\le\th\le1$ and
$\frac{1}{q}=\frac{1-\th}{p}+\frac{\th}{p'}$. Then for any map $u:
E\to H_q^c$ the following assertions are equivalent:
 \begin{enumerate}[{\rm i)}]
\item  $u$ is completely bounded.
\item  There is a constant $K$ such that for all finite
sequences $(a_k)\subset E$ and  $(\mu_k)\subset\real_+$
 \beq\label{LG2}
 \sum_k\|u(a_k)\|^2\le K^2 \Big[
 (1-\th)\|\sum_k\mu_k^{\th}a_k^*a_k\|_{p/2}+
 \th\,\|\sum_k\mu_k^{-(1-\th)}a_ka_k^*\|_{p/2}\Big].
 \eeq
\item  There are two positive unit elements $f, g\in
\big(L_{p/2}(M)\big)^*$ such that
 \beq\label{LG3}
 \|u(a)\|\le K\, \big(f(a^*a)\big)^{(1-\th)/2}
 \big(g(aa^*)\big)^{\th/2}\ ,\quad a\in E.
 \eeq
\end{enumerate}
Moreover, if $K$ denotes the best constant in (0.6) and (0.7),
then
 $$c_{p,\,\th}^{-1}\;K\le\|u\|_{cb}\le  K,$$
where $c_{p,\,\th}$ is a positive constant depending only on $p$
and $\th$, which can be controlled by an absolute constant.
\end{thm}

As in \cite{pisshlyak} and \cite{pis-Rp}, the key ingredient for
the proofs of Theorems~0.1 and 0.2 is some noncommutative
Khintchine type inequalities for Shlyakhtenko's generalized
circular systems. Note that the von Neumann algebra generated by a
generalized circular system is of type III. This forces us to work
with noncommutative $L_p$-spaces based on general von Neumann
algebras. These inequalities are presented in section~\ref{khin}.
The $L_\8$ case was already obtained in [PS]. We should emphasize
that the consideration of type III von Neumann algebras is
inevitable both for \cite{pisshlyak} and the present paper.

Besides these Khintchine type inequalities, we will still need two
tools. The first one is a new tensor product. This is defined in a
way  similar to the usual Haagerup tensor product, replacing the
row and column spaces $R$ and $C$ by their $L_p$-space
counterparts $R_p$ and $C_p$, the $p$-row and $p$-column spaces.
This new tensor product shares many properties with the usual
Haagerup tensor product. It seems notably interesting when
restricted to the subspaces of noncommutative $L_p$-spaces for
$p\ge 2$. This is developed in section~\ref{haag}, where the main
result is a characterization of maps factorable through $R_p$ (see
Theorem \ref{dual of hp}). This is the $L_p$-space analogue of the
well known Christensen-Sinclair's factorization for completely
bounded bilinear maps.

The second tool needed is the vector-valued noncommutative
$L_p$-space theory developed by Pisier \cite{pis-asterisque} for
injective semifinite von Neumann algebras, and especially, the
recent extension by Junge in \cite{ju-fubini} and \cite{ju-qwep}
to QWEP algebras. As said previously, the von Neumann algebra
generated by a generalized circular system is of type III. It is
non injective. However, it is QWEP. This explains why we really
need Junge's extension of Pisier's theory. Junge's work is briefly
discussed in section~\ref{vvLp}.

Sections~\ref{pfG} and \ref{pfLG} are devoted to the proofs of
Theorems~0.1 and 0.2, respectively. The last section contains some
applications. We mention here two of them. The first one is that
any completely bounded map from a subspace of a noncommutative
$L_p(M)$ into a quotient of $L_{p'}(M)$ has completely bounded
approximation property ($1/p+1/p'=1;\; 2<p<\8$). The second is a
characterization of (completely) bounded Schur multipliers from
$S_p$ (or a suitable subspace) to its dual ($2<p<\8$). The latter
is again the $L_p$-space analogue of the corresponding result in
\cite{pisshlyak}.


\section{Preliminaries}\label{prelim}

 In this section we collect some preliminaries necessary
to the whole paper. For clarity we divide the section into several
subsections.

\medskip\n {\bf 1.1.} {\it Operator spaces} :

\medskip  We will use standard notions and notation from operator
space theory. Our references are \cite{er-book} and
\cite{pis-intro}. $M_n$ will denote the algebra of all complex
$n\times n $ matrices. Given an operator space $E\subset B(H)$ we
denote by $M_n(E)$ the space of all $n\times n $ matrices with
entries in $E$. $M_n(E)$ is also an operator space equipped with
the operator space structure induced by that of $M_n(B(H))\cong
B(\ell_2^n(H))$. Let $u: E\to F$ be a linear map between two
operator spaces. $u$ is said to be {\it completely bounded} (c.b.
in short) if
 $$\|u\|_{cb}=\sup_n\|I_{M_n}\ot u\|_{M_n(E)\to M_n(F)}<\8.$$
Let $CB(E, F)$ denote the operator space of all c.b. maps from $E$
to $F$. $u$ is said to be {\it completely isomorphic} if $u$ is an
isomorphism and both $u, \; u^{-1}$ are c.b.. Similarly, we define
{\it complete contraction, complete isometry}.

Let $E, F, G$ be operator spaces. A bilinear map $u: E\times F \to
G$ is said to be {\it jointly completely bounded} (j.c.b. in
short) if the associated linear map $\wt u: E\to CB(F, G)$ is
c.b.. Then we set $\|u\|_{jcb}=\|\wt u\|_{cb}$. Let $JCB(E, F; G)$
denote the space of all j.c.b. maps from $E\times F$ to $G$. In
particular, a bilinear form $u: E\times F \to \comp$ is j.c.b. iff
the associated linear map $\wt u: E\to F^*$ is c.b..

Given an operator space $E$, we denote by $\overline{E}$ the
complex conjugate of $E$. As a vector space, $\overline{E}$ is the
same as $E$ but with the conjugate multiplication by a complex
scalar. If $x\in E$, $\bar x$ denotes the same vector $x$
considered as an element of $\overline{E}$. The norm of
$M_n(\overline{E})$ is defined by $\|(\bar
x_{ij})\|_{M_n(\overline{E})}=\|(x_{ij})\|_{M_n({E})}$. Thus the
map $x\mapsto \bar x$ establishes a complete anti-isometry between
$E$ and $\overline{E}$.

\medskip\n{\bf 1.2.} {\it Noncommutative $L_p$-spaces and
their natural operator space structure} :

\medskip
It is well known by now that there are several equivalent
constructions of noncommutative $L_p$-spaces associated with a von
Neumann algebra. In this paper we will use Haagerup's
noncommutative $L_p$-spaces (cf. \cite{haag-Lp}). \cite{terp-Lp}
is our main reference for these spaces. Throughout this paper $M$
will denote a general von Neumann algebra.

Let $M$ be a von Neumann algebra. For $0<p<\infty$, the spaces
$L_p(M)$ are constructed as spaces of measurable operators
relative not to $M$ but to a certain semifinite super von Neumann
algebra $\M$, namely, the crossed product of $M$ by one of its
modular automorphism groups. Let $(\theta_s)$ be the dual
automorphism group on $\M$. It is well known that $M$ is a von
Neumann subalgebra of $\M$ and that the position of $M$ in $\M$ is
determined by the group $(\theta_s)$ in the following sense:
 $$\forall x\in\M, \quad x\in M
 \Leftrightarrow(\forall s\in\real,\ \theta_s(x)=x).$$
Moreover, $\M$ is semifinite and can be canonically equipped with
a normal semifinite faithful trace $\tau$ such that
 $$\tau\circ \theta_s=e^{-s}\tau.$$

Let $L_0(\M,\tau)$ be the topological $\ast$-algebra of measurable
operators associated with $(\M,\tau)$ (in Nelson's sense
\cite{nelson}; see also \cite{terp-Lp}). The automorphisms
$\theta_s$, $s\in\real$,  extend to automorphisms of
$L_0(\M,\tau)$. For $0<p\le\infty$, the space $L_p(M)$ is defined
by
 $$L_p(M)=\{h\in L_0(\M,\tau)\mid \forall\;s\in\real\
 \theta_s(h)=e^{-s/p}h\}.$$
The space $L_\infty(M)$ coincides with $M$ (modulo the inclusions
$M\subset\M\subset L_0(\M,\tau)$). The spaces $L_p(M)$  are closed
self-adjoint linear subspaces of $L_0(\M,\tau)$. They are closed
under left and right multiplications by elements of $M$. If
$h=u|h|$ is the polar decomposition of $h\in L_0(\M,\tau)$, then
 $$h\in L_p(M)\Leftrightarrow u\in M \mbox{ and }|h|\in L_p(M).$$

It was shown by Haagerup that there is a linear homeomorphism
$\o\mapsto h_\o$ from $M_*$ onto $L_1(M)$ (equipped with the
vector space topology inherited from $L_0(\M,\tau) $), and this
homeomorphism preserves the additional structures (conjugation,
positivity, polar decomposition, action of $M$). It permits to
transfer the norm of $M_*$ into a norm on $L_1(M)$, denoted by
$\|\;\|_1$.

The space $L_1(M)$ is equipped with a distinguished bounded
positive linear functional $\tr$, the ``trace'', defined by
 $$\tr\,(h_\o)=\o(1),\quad \o\in M_*.$$
Consequently, $\|h\|_1=\tr\,(|h|)$ for every $h\in L_1(M)$.

For every $0<p<\infty$, the Mazur map $\M_+\to\M_+$, $x\mapsto
x^p$ extends by continuity to a map $L_0(\M,\tau)_+\to
L_0(\M,\tau)_+$, $h\mapsto h^p$ (cf. \cite{ray-ultra}). Then
 $$\forall\; h\in L_0(\M,\tau)_+\;,\quad
 h\in L_p(M)\Leftrightarrow h^p\in L_1(M).$$
For $h\in L_p(M)$ set $\|h\|_p=\|\,|h|^p\|_1^{1/p}$. Then
$\|\;\|_p$ is a norm or a $p$-norm according to  $1\leq p<\infty$,
or $0<p<1$. The associated vector space topology coincides with
that inherited from $L_0(\M,\tau)$.

Another important link between the spaces $L_p(M)$ is the {\it
external product}: in fact, the product of $L_0(\M,\tau)$,
$(h,k)\mapsto h\cdot k$, restricts to a bounded bilinear map
$L_p(M)\times L_q(M)\to L_r(M)$, where
$\frac{1}{r}=\frac{1}{p}+\frac{1}{q}$. This bilinear map has norm
one, which amounts to saying that the usual H\"older inequality
extends to  Haagerup $L_p$-spaces (called ``noncommutative
H\"older inequality'').

Assume that $\frac{1}{p}+\frac{1}{p'}=1$. Then the bilinear form
$L_p(M)\times L_{p'}(M)\to \comp$, $(h,k)\mapsto \tr\,(h\cdot k)$
defines a duality bracket between $L_p(M)$ and $L_{p'}(M)$, for
which $L_{p'}(M)$ is (isometrically) the dual of $L_p(M)$ (if
$p\neq\infty$); moreover we have the tracial property:
 $$\tr\,(hk)=\tr\,(kh),\quad h\in L_p(M), k\in L_{p'}(M).$$
In the sequel $\tr$ will denote this tracial functional associated
to any von Neumann algebra. In case of possible ambiguity, we will
write $\tr_M$ to indicate that the von Neumann algebra in
consideration is $M$.

\medskip

Now we turn to describe the natural operator space structure on
$L_p(M)$ as introduced in \cite{pis-asterisque} and
\cite{pis-intro}  (see also \cite{ju-fubini}). For $p=\infty$,
$L_\8(M)=M$ has its natural operator space structure as a  von
Neumann algebra. This yields an operator space structure on $M^*$,
the standard dual of $M$. Let us consider the case of $p=1$.
Recall that $L_1(M)$ coincides with the predual $M_*$ of $M$ at
the Banach space level. Thus one would attempt to define the
operator space structure on $L_1(M)$ as the one induced by that of
$M^*$ via the natural embedding $M_* \hookrightarrow
M^*=(M_*)^{**}$ (again as Banach spaces). However, as explained in
\cite [Chapiter 7] {pis-intro} , it is more convenient to consider
$L_1(M)$ as the predual of the opposite von Neumann algebra
$M^{op}$, which is isometric (but in general not completely
isomorphic) to $M$, and to equip $L_1(M)$ with the operator space
structure inherited from $(M^{op})^*$. One of the main reasons for
this choice is that it insures that the equality $L_1(M_n\otimes
M)=S_1^n\widehat \otimes L_1(M)$ (operator space projective tensor
product) holds true (see \cite{ju-fubini}). Finally, the operator
space structure of $L_p(M)$ is obtained by complex interpolation,
using the well known interpretation of $L_p(M)$ as the
interpolation space $(M,\;L_1(M))_{1/p}$ (see \cite{terp-int}). If
$M$ admits a normal faithful state $\f$, we can also use Kosaki's
interpolation \cite{ko-int}. Note that in this case all injections
$I_{\eta}$ of $M$ into $L_1(M)$ ($0\le\eta\le1$) considered in
\cite{ko-int} give completely isometric interpolation spaces,
exactly as in the Banach space level.

When $M$ is semifinite, we will always consider $L_p(M)$ as the
usual $L_p$-space constructed from a normal semifinite faithful
trace. We refer, for instance, to the survey \cite{px-survey} for
semifinite noncommutative $L_p$-spaces and for more  references.

\medskip\n{\bf 1.3.} {\it Vector-valued Schatten classes} :

\medskip  One of our main tools is the theory of vector-valued
noncommutative $L_p$-spaces. This theory was first introduced and
developed by  Pisier \cite{pis-asterisque} for injective
semifinite von Neumann algebras.  Pisier's theory can be easily
extended to general injective (so not necessarily semifinite) von
Neumann algebras. Very recently, Junge \cite{ju-fubini} and
\cite{ju-qwep} has partly extended this theory to QWEP von Neumann
algebras. We will discuss Junge's extension in some more details
later in section~\ref{vvLp}. Here we content ourselves only with a
brief description of vector-valued Schatten classes. The reader is
referred to \cite{pis-asterisque} for more information.

The Schatten classes $S_p$ are equipped with their natural
operator space structure as described previously. Now let $E$ be
an operator space. We define $S_1[E]$ as the operator space
projective tensor product $S_1 \widehat \otimes E$. Then for any
$1<p<\8$ we define $S_p[E]$ by interpolation:
 $$S_p[E]=\big(S_\8[E],\ S_1[E]\big)_{1/p}\ .$$
Note that in this interpolation formula, $S_\8[E]$ can be replaced
by $B(\ell_2)\ot_{\min}E$, namely, we have
 $$S_p[E]=\big(B(\ell_2)\ot_{\min}E,\ S_1[E]\big)_{1/p}\ .$$
By reiteration, for any $1\le p_0, p_1\le\8$ and $0<\th<1$
 $$\big(S_{p_0}[E],\ S_{p_1}[E]\big)_{\th}=S_p[E]\ ,$$
where $\frac{1}{p}=\frac{1-\th}{p_0}+\frac{\th}{p_1}$. More
generally, given a compatible couple $(E_0, E_1)$ of operator
spaces we have (completely isometrically)
 \beq\label{interpolation Sp}
 \big(S_{p_0}[E_0],\ S_{p_1}[E_1]\big)_{\th}
 =S_p[(E_0,\ E_1)_\th]\ .
 \eeq

The usual duality for Schatten classes extends to the
vector-valued case too. Let $1\le p<\8$. Then
 \beq\label{dual Sp}
 \big(S_p[E]\big)^*=S_{p'}[E^*]\quad
 \mbox{completely isometrically}.
 \eeq
Here and throughout this paper, $p'$ always denotes the conjugate
index of $p$. The duality is given as follows. Let $x=(x_{ij})\in
S_p[E]$ and $\xi=(\xi_{ij})\in S_{p'}[E^*]$. Then
 $$\la \xi,\; x\ra=\Tr(\, ^t\xi\, x)=\sum_{ij}\xi_{ij}(x_{ij})\ ,$$
where $\Tr$ is the usual trace on $B(\ell_2)$. We call the
reader's attention to that in the theory of operator spaces the
duality between $S_p$ and $S_{p'}$ (in the scalar case) is given
by
 $$\la y,\; x\ra=\Tr(\, ^ty\, x)=\sum_{ij}x_{ij}y_{ij}\ .$$
This is consistent with the natural operator space structure on
$S_1$ described previously, which is the predual of
$B(\ell_2)^{op}$.

There is a nice description of the norm of $S_p[E]$ in terms of
that of $S_\8[E]$, given by \cite [Theorem 1.5] {pis-asterisque}:
for any $x=(x_{ij})\in S_p[E]$ we have
 \beq\label{Sp via S infinite}
 \|x\|_{S_p[E]}=\inf\big\{\|\a\|_{2p}\,\|y\|_{S_\8[E]}\,
 \|\b\|_{2p}\big\},
 \eeq
where the infimum runs over all factorizations $x=\a y\b$ with
$\a, \b\in S_{2p}$ and $y\in S_\8[E]$. Conversely, the norm of
$S_\8[E]$ can be recovered from that of $S_p[E]$ as follows (see
\cite[Lemma 1.7]{pis-asterisque}): for any $x=(x_{ij})\in S_\8[E]$
 \beq\label{S infinite via Sp}
 \|x\|_{S_\8[E]}=\sup\big\{\|\a y\b\|_{S_p[E]}\ :\ \a, \b
 \in S_{2p},\ \|\a\|_{2p}\le 1, \|\b\|_{2p}\le 1\big\}.
 \eeq
From (\ref{Sp via S infinite}) and (\ref{S infinite via Sp}) we
can easily deduce the following more general formula. Let $1\le
p<q\le\8$ and $\frac{1}{r}=\frac{1}{p}-\frac{1}{q}$. Then for any
$x=(x_{ij})\in S_q[E]$
  \beq\label{Sq via Sp}
  \|x\|_{S_q[E]}=\sup\big\{\|\a y\b\|_{S_p[E]}\ :\ \a, \b
 \in S_{2r},\ \|\a\|_{2r}\le 1, \|\b\|_{2r}\le 1\big\}.
 \eeq
(\ref{S infinite via Sp}) also implies the following convenient
characterization of c.b. maps, which will be repeatedly used later
(cf. \cite[Lemma 1.7]{pis-asterisque}).

\begin{lem}\label{cb via Sp}
 Let $E$ and $F$ be two operator spaces. Then a linear map
$u: E\to F$ is c.b. iff
 \[
 \sup_n\big\|I_{S_p^n}\ot u: S_p^n[E]\to
 S_p^n[F]\big\|<\infty;
 \]
moreover in this case the supremum above is equal to $\|u\|_{cb}$.
\end{lem}

More generally, if $H$ is a Hilbert space, we can analogously
define $S_p[H; E]$, the $E$-valued Schatten classes based on $H$.
$S_p[E]$ (resp. $S_p^n[E]$) corresponds to the case of infinite
dimensional separable $H$ (resp. $\dim H=n$). All preceding
properties still hold for general $S_p[H; E]$.

\medskip

We now specialize the discussion above to the case when $E$ is a
subspace  of a noncommutative $L_p(M)$. In this case, the theory
becomes much simpler. Note that there is a natural algebraic
identification of $L_p(M_n\otimes M)$ with $M_n(L_p(M))$. Then
$S_p^n[L_p(M)]$ is nothing but the linear space $M_n(L_p(M))$
equipped with the norm of $L_p(M_n\otimes M)$. More generally, if
$E\subset L_p(M)$ is a closed subspace, the norm on $S_p^n[E]$ is
induced by that of $S_p^n[L_p(M)]$. In the infinite dimensional
case, $S_p[L_p(M)]$ is completely isometrically identified with
$L_p(B(\ell_2)\overline{\otimes} M)$ for all $1\le p<\8$. If
$E\subset L_p(M)$, then $S_p[E]$ is the closure in
$L_p(B(\ell_2)\overline{\otimes} M)$ of the algebraic tensor
product $S_p\ot E$.

\medskip\n{\bf 1.4.} {\it Column and row spaces} :

\medskip  The column and row spaces, $C$ and $R$, play an
important role in the whole theory of operator spaces. Note that
$C$ and $R$ are respectively the (first) column and row subspaces
of $S_\8$. The $L_p$-space counterparts of $R$ and $C$ will play
an essential  role in the present paper. Let $C_p$ (resp. $R_p$)
denote the subspace of $S_p$ consisting of matrices whose all
entries but those in the first column (resp. row) vanish. So
$C_\8$ and $R_\8$ are just $C$ and $R$ respectively. It is clear
that $C_p$ and $R_p$ are completely 1-complemented subspaces of
$S_p$. We have the following completely isometric identities: for
any $1\le p\le\8$
 \beq\label{dual Cp-Rp}
 (C_p)^*\cong C_{p'}\cong R_p\quad\mbox{and}\quad
 (R_p)^*\cong R_{p'}\cong C_p.
 \eeq
$C_p$ and $R_p$ can  be also defined via interpolation from $C$
and $R$. We view $(C, R)$ as a compatible couple by identifying
both of them with $\ell_2$ (in the Banach space level!), i.e. by
identifying the canonical bases $(e_{k,1})$ of $C_p$ and
$(e_{1,k})$ of $R_p$ with $(e_k)$ of $\ell_2$. Then
 $$C_p=(C,\ R)_{1/p}=(C_\8,\ C_1)_{1/p}\quad\mbox{and}\quad
 R_p=(R,\ C)_{1/p}=(R_\8,\ R_1)_{1/p}\ .$$
By reiteration, for any $1\le p_0, p_1\le\8$ and $0<\th<1$
 \beq\label{interpolation Cp-Rp}
 C_p=(C_{p_0},\ C_{p_1})_{\th}\quad\mbox{and}\quad
 R_p=(R_{p_0},\ R_{p_1})_{\th}\; ,
 \eeq
where $\frac{1}{p}=\frac{1-\th}{p_0}+\frac{\th}{p_1}$. Like $C$
and $R$, $C_p$ and $R_p$ are also 1-homogenous 1-Hilbertian
operator spaces. We refer to \cite{pis-asterisque} for the proofs
of all these elementary facts.

\medskip

More generally, given a Hilbert space $H$ and $1\le p\le\8$ we
denote by $H_p^c$ (resp. $H_p^r$) the Schatten $p$-class
$S_p(\comp, H)$ (resp. $S_p(\bar H, \comp)$) equipped with its
natural operator space structure. When $H$ is separable and
infinite dimensional, $H_p^c$ and $H_p^r$ are respectively $C_p$
and $R_p$ above. If $\dim H=n<\8$, we set $H_p^c=C_p^n$ and
$H_p^r=R_p^n$. All properties for $C_p$ and $R_p$ mentioned above
hold for $H_p^c$ and $H_p^r$ too. We will call $H_p^c$ (resp.
$H_p^r$) a {\em $p$-column} (resp. {\em $p$-row}) {\em space}.

\medskip

Now let $E$ be an operator space. We denote by $C_p[E]$ (resp.
$R_p[E]$) the closure of $C_p\ot E$ (resp. $R_p\ot E$) in
$S_p[E]$. Again $C_p[E]$ and  $R_p[E]$ are completely
1-complemented subspaces of $S_p[E]$. If $E$ is a subspace of a
noncommutative $L_p(M)$, the norm of $C_p[E]$ is easy to be
determined. For any finite sequence $(x_k)\subset E$
 \beq\label{col Lp}
 \big\|\sum_kx_k\ot e_k\big\|_{C_p[E]}
 =\big\|\big(\sum_kx_k^*x_k\big)^{1/2}\big\|_{L_p(M)}\; ,
 \eeq
where $(e_k)$ denotes the canonical basis of $C_p$. More
generally, if $a_k\in C_p$, then
$$\big\|\sum_kx_k\ot a_k\big\|_{C_p[E]}
 =\big\|\big(\sum_k\la a_k, a_j\ra x_k^*x_j\big)^{1/2}\big\|_{L_p(M)}\; ,$$
where $\la\; ,\;\ra$ denotes the scalar product in $C_p$. (In
terms of matrix product, $\la a_k, a_j\ra=a_k^*a_j$.) We also have
a similar description for $R_p[E]$.

\medskip\n{\bf 1.5.} {\it Factorization through $C_p$ and $R_p$} :

\medskip  The following definition goes back to \cite{pis-OH}
(see also \cite{pis-intro}). Given operator spaces $X, E$ and $F$,
we denote by $\G_X(E, F)$ the family of all maps $u: E\to F$ which
factors through $X$, namely, all $u$ which admit  a factorization
$\displaystyle E\mathop{\longrightarrow}^\alpha X
\mathop{\longrightarrow}^\beta F$ with c.b. maps $\a$ and $\b$.
 For $u\in \G_X(E, F)$ define
 $$\g_{X}(u)=\inf\big\{\|\a\|_{cb}\|\b\|_{cb}\;:\; u=\b\a,\;
 \ \a\in CB(E, X),\; \b\in CB(X, F)\big\}.$$
We will need this  notably when $X$ is $R_p, C_p$ or $R_p\op_p
C_p$. Then $X$ can be any one of these spaces associated with an
arbitrary Hilbert space. Let us make this more precise. A map $u:
E\to F$  is said to {\it be factorable through $C_p$} if $u$
admits a factorization $E\displaystyle{\mathop{\longrightarrow}^\a
H_p^c\mathop{\longrightarrow}^\beta F}$ for some Hilbert space $H$
and c.b. maps $\a$ and $\b$. Let $\G_{C_p}(E, F)$ denote the space
of all maps between $E$ and $F$ factorable through $C_p$. For
$u\in \G_{C_p}(E, F)$ define
 $$\g_{C_p}(u)=\inf\big\{\|\a\|_{cb}\|\b\|_{cb}\big\},$$
where the infimum runs over all factorizations $u=\b\a$ as above.
Then $\big(\G_{C_p}(E, F),\g_{C_p}\big)$ is a Banach space.
Similarly, we define the factorization through $R_p$ and $C_p\op_p
R_p$, respectively. The resulting spaces are denoted respectively
by $\G_{R_p}(E, F)$ and $\G_{C_p\op_p R_p}(E, F)$ equipped with
$\g_{R_p}$ and $\g_{C_p\op_p R_p}$.


\section{ A Haagerup type tensor norm and
factorization through $p$-row spaces}\label{haag}

In this section we introduce a Haagerup type tensor norm. We refer
to \cite{er-book} and \cite{pis-intro} for the usual Haagerup
tensor product. This new tensor norm enjoys many properties of the
usual Haagerup tensor norm. It turns out especially satisfactory
when restricted to subspaces of noncommutative $L_p$-spaces for
$p\ge 2$. In this latter case we obtain a description of the
linear functionals continuous with respect to this tensor norm.
This is the $L_p$-space analogue of the well-known factorization
theorem due to Christensen and Sinclair.

Let us begin with some conventions. All row and column matrices
below contain only finitely many non zero coefficients, so they
can be considered as finite row and column matrices. When $a$ and
$ b$ are two row (resp. column) matrices with entries in $E$, $(a,
b)$ (resp. $ ^t(a, b)$) is again a row (resp. column) matrix with
coefficients in $E$.

\medskip

The definition below is a generalization of the Haagerup tensor
product, just replacing $R_\8$ and $C_\8$ in the usual Haagerup
tensor product by $R_p$ and $C_p$. Recall that if $a\in R_p[E]$ is
a row matrix and $b\in C_q[F]$ a column one, $a\odot b$ denotes
the element in $E\otimes F$ given by
 $$a\odot b=\sum_k a_k\otimes b_k.$$

\begin{defi}\label{def hpq} Let $E$ and $F$ be
operator spaces and $1\le p,q\le\8$. Given $x\in E\otimes F$ we
define
 $$\|x\|_{h_{p,q}}=\inf\big\{\|a\|_{R_p[E]}\,\|b\|_{C_q[F]}\ :\
 x=a\odot b,\ a\in R_p[E],\ b\in C_q[F]\big\}.$$
If $p=q$, $\|x\|_{h_{p,q}}$ is simply denoted by $\|x\|_{h_{p}}$.
\end{defi}

In general, $\|\,\|_{h_{p,q}}$ is not a norm, but only a
quasi-norm (see the proposition below). However, when $E$ and $F$
verify a certain 2-convexity, then $\|\,\|_{h_{p,q}}$ is indeed a
norm.

\begin{defi}\label{2-convex}
 Let $E$ be an operator space and $1\le
p\le\8$. $E$ is said to be $R_p$-$2$-convex $($resp.
$C_p$-$2$-convex$)$ if for any $a, b \in R_p[E]$ $($resp.
$C_p[E])$
 $$\|(a, b)\|_{R_p[E]}\le\big(\|a\|_{R_p[E]}^2
 + \|b\|_{R_p[E]}^2\big)^{1/2}$$
$($resp.
 $$\|\, ^t(a, b)\|_{C_p[E]}\le\big(\|a\|_{C_p[E]}^2
 + \|b\|_{C_p[E]}^2\big)^{1/2}\ )\ .$$
 \end{defi}

It is easy to see that  the family of $R_p$-2-convex (resp.
$C_p$-2-convex) operator spaces is closed with respect to
subspaces and quotients . It is also easy to check (using
(\ref{col Lp})) that subspaces of non-commutative $L_p$-spaces
($p\ge 2$) are $R_p$-2-convex and $C_p$-2-convex. Any operator
space is $R_\8$-2-convex and $C_\8$-2-convex.

\medskip

\n{\bf Remark.} Let $E$ (resp. $F$) be an $R_p$-2-convex (resp.
$C_q$-2-convex) operator space. Then the functional $\|\;
\|_{h_{p,q}}$ defined by Definition \ref{def hpq} is a $\g$-norm
in the sense of \cite{pis-OH}. Indeed, given a positive element
$u=\sum a_k\ot\bar a_k$ in $E\ot\overline E$, define
 $$r_p(u)=\big\|(a_1, a_2,\cdots )\big\|_{R_p[E]}^2\;.$$
Then the $R_p$-2-convexity of $E$ guarantees that $r_p$ is a
weight on $(E\ot\overline E)_+$ in Pisier's sense. Similarly, we
have a weight $c_q$ on $(F\ot\overline F)_+$ corresponding to the
norm of $C_q[F]$. Then $\|\; \|_{h_{p,q}}$ is exactly the
$\g$-norm associated to $r_p$ and $c_q$ defined in
\cite[section~6] {pis-OH}.

\medskip

The following is elementary. As usual, an element $x\in E\ot F$ is
also regarded as a map from $E^*$ to $F$. Its adjoint $^tx$ is
from $F^*$ to $E$. Then the norm $\|\, ^tx\|_{h_{p,q}}$ is the
norm of $x$ in $F\otimes_{h_{p,q}} E$.

\begin{prop}\label{hpq}
Let $E$ and $F$ be operator spaces and $1\le p,q\le\8$.
 \begin{enumerate}[{\rm i)}]
 \item  $\|\ \|_{h_{p,q}}$ is a quasi-norm on $ E\otimes F$. If in
addition $E$ and $F$ are $R_p$-2-convex and $C_q$-2-convex,
respectively, then $\|\ \|_{h_{p,q}}$ is a norm.
 \item  For any $x\in E\otimes F$ there are $a=(a_1, ... ,
 a_n)\in R_p[E]$ and $b=~ ^t(b_1, ... ,
 b_n)\in C_q[F]$ such that  both $(a_1, ... ,
 a_n)$ and $(b_1, ... , b_n)$ are linearly independent and such that
  $$\|x\|_{h_{p,q}}=\|a\|_{R_p[E]}\;\|b\|_{C_q[F]}\ .$$
Moreover, $a$ and $b$ can be chosen to further satisfy
 $$\|\, ^tx\|_{h_{p,q}}=\|\, ^tb\ \Delta^{-1}\|_{R_p[F]}\,
 \|\Delta\, ^ta\|_{C_q[E]}\;,$$
where $\Delta$ is a  positive definite diagonal matrix.
 \end{enumerate}
 \end{prop}

\pf  i) We have the following quasi-triangle inequality:
 $$\|x+y\|_{h_{p,q}}\le 2\big(\|x\|_{h_{p,q}}
 +\|y\|_{h_{p,q}}\big).$$
This is proved by a standard argument that is left to the reader.
The 2-convexity assumption implies  the validity of the triangle
inequality. Thus it remains to check that $\|x\|_{h_{p,q}}=0$
implies $x=0$. To this end we show $\|x\|_{\e}\le
\|x\|_{h_{p,q}}$, where $\|\,\|_{\e}$ denotes the Banach space
injective tensor norm. Let $\xi\in E^*$ and $\eta \in F^*$ be unit
vectors. Then for any factorization $x=a\odot b$
 $$|\la \xi\ot\eta,\ x\ra|=\big|\sum_k \xi(a_k)\,\eta(b_k)\big|
 \le \Big(\sum_k|\xi(a_k)|^2\Big)^{1/2}\;
 \Big(\sum_k|\eta(b_k)|^2\Big)^{1/2}\;.$$
Let $(\a_1, \a_2, \cdots)\in\el_2$ be a unit vector such that
 $$\Big(\sum_k|\xi(a_k)|^2\Big)^{1/2}=\sum_k\a_k\xi(a_k)
 \le \big\|\sum_k\a_k a_k\big\|.$$
Writing
 $$\sum_k\a_k a_k= (a_1, a_2, \cdots)\; ^t(\a_1, \a_2, \cdots),$$
we get
 $$\big\|\sum_k\a_k a_k\big\|\le \|a\|_{R_p[E]}\;
 \big\|\,^t(\a_1, \a_2, \cdots)\big\|_{B(\el_2)}\le \|a\|_{R_p[E]}\;.$$
Thus
 $$\Big(\sum_k|\xi(a_k)|^2\Big)^{1/2}\le \|a\|_{R_p[E]}\;.$$
Similarly,
 $$\Big(\sum_k|\eta(b_k)|^2\Big)^{1/2}\le \|b\|_{C_q[E]}\;.$$
Hence it follows that $\|x\|_{\e}\le \|x\|_{h_{p,q}}$.

ii) Recall that the $R_p$- and $C_q$-norms satisfy the following
elementary property: for any scalar matrix $\a\in B(\ell_2)$
 $$\|a\a\|_{R_p[E]}\le \|a\|_{R_p[E]}\,\|\a\|\quad\mbox{and}\quad
 \|\a b\|_{C_q[F]}\le \|\a\|\, \|b\|_{C_q[F]}\ .$$
Using this property one can easily prove the first assertion of
ii) exactly as in the case of the usual Haagerup tensor product.
See \cite[section~9.2] {er-book} for more details. Similarly, the
second part can be proved as \cite[Proposition 1.7]{pisshlyak}.
\cqd

\medskip

We denote by $E\otimes_{h_{p,q}} F$ the completion of ($E\otimes
F$, $\|\; \|_{h_{p,q}})$. Again $E\otimes_{h_{p,q}} F$ is denoted
by $E\otimes_{h_{p}} F$ in the case of $p=q$. By definition, the
tensor product $E\otimes_{h_{p,q}} F$ is projective. Proposition
\ref{hpq} ii) implies that it is also injective.

\medskip

From now on we consider only the case where $p=q\ge 2$, and
specialize the above tensor product to subspaces of noncommutative
$L_p$-spaces. (Recall that all these subspaces are both
$R_p$-2-convex and $C_p$-2-convex.) For these spaces we will have
a satisfactory description of the dual space of $E\otimes_{h_p}F$.
We first need to characterize the c.b. maps from a subspace of a
noncommutative $L_p$ to $R_p$. The following is the $L_p$-space
analogue of a well-known result on maps with values in $R$ due to
Effros and Ruan (cf. \cite{er-self}). The main point here is the
implication iv) $\Rightarrow$ i). We should emphasize that this
result (without the assertion ii) below) is a special case of
Theorem \ref{LG}, corresponding to the endpoint cases $\th=0$ and
$\th=1$. In this special case, the arguments are elementary and
much simpler than that for the general case as in Theorem
\ref{LG}.

\begin{prop}\label{to Rp}
 Let $E\subset L_p(M)$ be a closed
subspace $(2< p<\8)$ and $H$ a Hilbert space. Let $u: E\to H_p^r$
be a linear map. Then the following assertions are equivalent
\begin{enumerate}[{\rm i)}]
 \item  $u$ is c.b.
 \item  $I_{R_p}\otimes u$ extends to a bounded map from
$R_p[E]$ to $ R_p[H_p^r]$.
\item  There is a constant $c$ such that for any finite
sequence $(a_k)\subset E$
 \beq\label{to Rp3}
 \Big(\sum_k\|u(a_k)\|^2\Big)^{1/2} \le c \,
 \big\|\big(\sum_ka_ka_k^*\big)^{1/2}\big\|_{p}\;.
 \eeq
\item  There is  a positive unit element $f\in
\big(L^{p/2}(M)\big)^*=L_{(p/2)'}(M)$ such that
 \beq\label{to Rp4}
 \|u(a)\|\le c\,\big(f(aa^*)\big)^{1/2}\ ,\quad a\in E.
 \eeq
 \end{enumerate}
 Moreover, if one of these assertions holds, $\|u\|_{cb}$,
$\|I_{R_p}\otimes u\|$ and the smallest constants $c$ in (\ref{to
Rp3}) and (\ref{to Rp4}) are all equal, and $u$ admits a c.b.
extension to $L_p(M)$ with the same c.b. norm.

 \n We have a similar result for maps with values in $H_p^c$ $($with
necessary changes in {\rm ii)-iv)} above$)$.
\end{prop}

\pf i) $\Rightarrow$ ii). By Lemma 1.1, $u$ is c.b. iff
$I_{S_p}\otimes u$ extends to a bounded map from $S_p[E]$ to
$S_p[H_p^r]$. In particular, i) implies ii).

 ii) $\Rightarrow$ iii).  Assume ii). Then for any finite row
$a=(a_1, a_2,...)\in R_p[E]$
 $$\|I_{R_p}\otimes u(a)\|_{R_p[H_p^r]}\le c\
 \|a\|_{R_p[E]}\ .$$
However,
 $$\|I_{R_p}\otimes u(a)\|_{R_p[H_p^r]}
 =\big(\sum_k\|u(a_k)\|^2\big)^{1/2}
 \quad\mbox{and}\quad
 \|a\|_{R_p[E]}=\|(\sum_ka_ka_k^*)^{1/2}\|_p\ .$$
Thus (\ref{to Rp3}) follows.

iii) $\Rightarrow$ iv). This can be done by a standard application
of the Hahn-Banach theorem. We give a sketch of the proof for
completeness. Let $S$ denote the positive cone of the unit ball of
$\big(L_{p/2}(M)\big)^*=L_{(p/2)'}(M)$. Then $S$ is a compact
space when equipped with the w*-topology. Given a finite sequence
$(a_k)\subset E$ set
 $$\f_{(a_k)}(s)=c^2\,s(\sum_ka_ka_k^*)-\sum_k\|u(a_k)\|^2\;,\quad
 s\in S.$$
Then $\f_{(a_k)}$ is a continuous function on $S$ and (\ref{to
Rp3}) implies that its maximum is positive. Let $A$ denote the
closure of the family  of all such functions $\f_{(a_k)}$. Then
$A$ is a closed convex cone and disjoint from $A_-$, where $A_-$
is the open convex cone of all negative continuous functions on
$S$. Therefore, by the Hahn-Banach theorem there is a probability
measure $\mu$ on $S$ such that
 $$\int_S\big(c^2s(aa^*) - \|u(a)\|^2\big)\, d\mu(s)\ge 0.$$
Define
 $$f(a)=\int_S s(a) \, d\mu(s),\quad a\in L_{p/2}(M).$$
Then $f$ is a positive unit element in  $L_{(p/2)'}(M)$ satisfying
(\ref{to Rp4}).

iv) $\Rightarrow$ i). Note that $\la b, a\ra=f(ab^*)$ defines a
semi-scalar product on $E$. Quotiented by its kernel, $E$ becomes
a pre-Hilbert space whose completion is denoted by $K$.  It is
clear that the identity on $E$ induces a contractive inclusion of
$E$ into $K$, denoted by $i_E$. On the other hand, (\ref{to Rp4})
implies that there is a bounded operator $\wh u: K\to H$ such that
$\|\wh u\|\le c$ and $u=\wh u\circ i_E$. We now equip $K$  with
the operator space structure of $K_p^r$, i.e. we consider $K$  as
$K_p^r$ in the category of operator spaces.  Then by the
homogeneity of $p$-row spaces,  $\wh u: K_p^r\to H_p^r$ is
automatically c.b. and $\|\wh u\|_{cb}=\|\wh u\|$. Therefore, it
remains to show that $i_E: E\to K_p^r$ is completely contractive.
By Lemma 1.1 it then suffices to prove
 $$\|I_{S_p^n}\otimes i_E\ :\ S_p^n[E]\to S_p^n[K_p^r]\,\|\le
 1,\quad \forall n\ge 1.$$
To this end we first need to identify the action of an operator in
$S_p^n[K_p^r]$.  Recall that $K_p^r=S_p(\overline{K}, \comp)$, and
so $S_p^n[K_p^r]=S_p(\ell_2^n(\overline{K}), \ell_2^n)$. Let
$y=(y_{i j})\in S_p^n[K_p^r]$. Considered as an operator in
$S_p(\ell_2^n(\overline{K}), \ell_2^n)$,   $y$ and $y^*$ act as
follows: for any $\b=(\b_k)\in \ell_2^n(\overline{K})$ and
$\a=(\a_k)\in \ell_2^n$
 $$y(\b)=\big(\sum_j\la\b_j,\ y_{i j}\ra\big)_{1\le i\le n}
 \quad\mbox{and}\quad
 y^*(\a)=\big(\sum_i \a_i\overline{y}_{i j}
 \big)_{1\le j\le n}\ .$$
Then it is easy to deduce that the $n\times n$ complex matrix $y
y^*$ is given by
 $$yy^*=\big(\sum_k\ \la y_{j k},\ y_{i k}\ra\big)_{1\le i\,j\le n}\ .$$
(Here $\la \cdot,\ \cdot\ra$ stands for the scalar product in
$K$.)

Now let $x=(x_{ij})\in S_p^n[E]$ and $y=I_{S_p^n}\otimes i_E(x)$.
Then the discussion above yields
 $$yy^*=\big(\sum_k\ f( x_{i k}x_{j k}^*)\ \big)_{1\le i,j\le n}\ .$$
Let $\a\in S_{(p/2)'}^n$ of norm 1. Then
 $$\Tr(\a yy^*)=\Tr\otimes \tr\big((\a\otimes f)(xx^*)\big)
 \le\|\a\otimes f\|_{(p/2)'} \|xx^*\|_{p/2}\le\|xx^*\|_{p/2}\ .$$
Taking the supremum over all such $\a$, we obtain
 $$\|yy^*\|_{p/2}\le \|xx^*\|_{p/2}\ ,$$
and so
 $$\|I_{S_p^n}\otimes i_E(x)\|_{S_p^n[K_p^r]}
 \le \|x\|_{S_p^n[E]}\ .$$
Therefore, $i_E$ is completely contractive. Thus we have proved
the equivalence between i)-- iv).

Finally, suppose one of i)-- iv) holds. Then from the previous
arguments we see that all the relevant constants are equal.
Moreover, from the proof of iv) $\Rightarrow$ i), $u$ factors as
$\displaystyle E\mathop{\longrightarrow}^{i_E}
K_p^r\mathop{\longrightarrow}^{\wh u}H_p^r$. From this one easily
deduces that $u$ admits a c.b. extension to the whole $L_p(M)$.
Indeed, let $\wt K$ be the Hilbert space constructed from $L_p(M)$
with respect to the semi-scalar product $\la b, a\ra=f(ab^*)$ as
previously.  Then $\wt K$ contains $K$ as an isometric subspace.
Let $i_{L_p(M)}: L_p(M)\to {\wt K}_p^r$ be the natural inclusion.
Then $i_{L_p(M)}$ is completely contractive. Let $P_K: \wt K\to K$
be the orthogonal projection. By the homogeneity of ${\wt K}_p^r$,
$P_K$ is completely contractive. Then $U=\wh u\, P_K\,i_{L_p(M)}$
is the desired extension of $u$.\cqd

\medskip

\n{\bf Remarks.} i) The equivalence in Proposition \ref{to Rp}
does not hold for general $R_p$-2-convex spaces. In fact,  a
simple interpolation argument shows that if $2\le q< p$, any
noncommutative $L_q$ is $R_p$-2-convex. Consequently, $R_q$ is
$R_p$-2-convex for such $p, q$. However, one can prove that the
complete boundedness of a map $u: R_q\to R_p$ is not equivalent to
the boundedness of $I_{R_p}\ot u: R_p[R_q]\to R_p[R_p]$. On the
other hand, $R_q$ is a quotient of a subspace of $R_p\op_p C_p$
(see \cite{xu-emb}). Thus this example also shows that Proposition
\ref{to Rp} does not hold for quotients of subspaces of
noncommutative $L_p$-spaces too.

ii) Using the weight $r_p$ introduced in the remark following
Definition \ref{2-convex}, part iii) of Proposition \ref{to Rp}
can be rephrased as  $\pi_{2,\,r_p}(u)\le c$, where
$\pi_{2,\,r_p}$ is the $(2, r_p)$-summing norm defined in
\cite[section 5] {pis-OH}. Thus for any $u$ as in Proposition
\ref{to Rp}, we have $\|u\|_{cb}=\pi_{2,\,r_p}(u)$.

\medskip

The proof of the implication iii) $\Rightarrow$ i) of Theorem
\ref{LG} in section~5 gives an alternative proof of the
implication iv) $\Rightarrow$ i) above (corresponding to $\th=1$
in Theorem \ref{LG}). The proof given above has an advantage that
it yields the natural factorization of $u$ by a kind of change of
density as in the commutative case. The natural inclusion $i_E:
E\to K_p^r$ constructed previously  will be used several times in
the sequel. For later reference let us explicitly record this as
follows.

\begin{rk}\label{iE}
 Let $E\subset L_p(M)$ be a closed subspace with $p>2$ and $f$ a
positive unit functional on $L_{p/2}(M)$. Let $K$ be the Hilbert
space obtained from $E$ relative to the semi-scalar product $\la
b,\; a\ra=f(ab^*)$. Then the natural inclusion $i_E: E\to K_p^r$
is completely contractive. A similar statement holds for the
semi-scalar product $\la b,\; a\ra=f(b^*a)$ (the resulting
Hilbertian operator space is then a $p$-column space).
\end{rk}

\begin{thm}\label{dual of hp}
 Let $M$ be a von Neumann algebra  and
$E, F\subset L_p(M)$ be closed subspaces $(2< p<\8)$. Let $u:
E\otimes F\to\comp$ be a linear functional and $c>0$ a constant.
The following assertions are equivalent
 \begin{enumerate}[{\rm i)}]
 \item  $u$ defines  a continuous functional on
$E\otimes_{h_p}F$ of norm $\le c$.
 \item  For all finite sequences $(a_k)\subset E$ and $(b_k)\subset F$
 \beq\label{dual of hp2}
 \big|\sum_k u(a_k\otimes b_k)\big|\le c\,
 \big\|\big(\sum_ka_ka_k^*\big)^{1/2}\big\|_p\,
 \big\|\big(\sum_kb_k^*b_k\big)^{1/2}\big\|_p\;.
 \eeq
 \item  There are positive unit functionals
$f, g\in \big(L_{p/2}(M)\big)^*$ such that
  \beq\label{dual of hp3}
  |u(a\otimes b)|\le c\, \big(f(aa^*)\, g(b^*b)\big)^{1/2}\ , \quad a\in
 E,\ b\in F.
 \eeq
 \item  The associated linear
map $\widetilde u: E\to F^*$ belongs to $\G_{R_p}(E, F^*)$ and
$\g_{R_p}(\widetilde u)\le c$.
\end{enumerate}
Moreover, if one of these assertions holds, $u$ has an extension
to $L_p(M)\otimes_{h_p} L_p(M)$ with the same norm.

\n We have a similar result for maps belonging to $\G_{C_p}(E,
 F^*)$.
\end{thm}

\pf  Going back to the definition of the norm of
$E\otimes_{h_p}F$, we see that ii) is just a reformulation of i).
The implication ii) $\Rightarrow$ iii) is shown by a standard
argument using the Hahn-Banach theorem as in the proof of iii)
$\Rightarrow$ iv) of Proposition \ref{to Rp}. Conversely, iii)
$\Rightarrow$ ii) is a simple consequence of the H\"older
inequality. It remains to show the equivalence iii)
$\Leftrightarrow$ iv).

First assume iii). Let $H$ be the Hilbert space obtained from $E$
relative the  semi-scalar product  $\la a,\; a'\ra=f(a'a^*)$ (see
Remark \ref{iE}). Similarly, let $K$ be the Hilbert space
associated with $F$ and the semi-scalar product $\la b,\; b'\ra=
g(b^*b')$. Let $i_E$ and $i_F$ be the natural inclusions of $E$
into $H$, respectively, of $F$ into $K$. Then (\ref{dual of hp3})
implies that there is a bounded operator $\wh u: H\to
\overline{K}$ with $\|\wh u\|\le c$ such that
 $$u(a, b)=\la\,\overline{\wh u\,i_E(a)} ,\ i_F(b)\,\ra,
 \quad a\in E,\ b\in F.$$
Thus we deduce
 $$\widetilde u=i_F^*\,\wh u\,i_E.$$
We now equip $H$ (resp. $K$)  with the operator space structure of
$H_p^r$ (resp. $K_p^c$). Then by Remark \ref{iE}, $i_E: E\to
H_p^r$ and $i_F: F\to K_p^c$ are completely contractive. Thus
$i_F^*: {\overline K}_p^r\to F^*$ is also completely contractive.
On the other hand,  $\wh u: H_p^r\to {\overline K}_p^r$ is c.b.
and has $\|\wh u\|$ as its cb-norm. Set $\a=i_E$ and
$\b=i_F^*\,\wh u$. Then $\a\in CB(E,H_p^r),\ \b\in CB(H_p^r,
F^*)$, $\widetilde u=\b\a$ and $\|\a\|_{cb}\|\b\|_{cb}\le c$. Thus
$\widetilde u\in \G_{R_p}(E, F^*)$ and $\g_{R_p}(\widetilde u)\le
c$.

Conversely, assume iv). By Proposition \ref{to Rp}, it is not hard
to see that $u$ satisfies ii). Therefore, the equivalence between
i) -- iv) has been proved.

Moreover, if iv) is verified, then by Proposition \ref{to Rp}, $u$
extends to $L_p(M)\otimes_{h_p} L_p(N)$ with the same norm. \cqd

\medskip

Theorem \ref{dual of hp} has the following extension with almost
the same proof.

\begin{rk}\label{dual of hpq}
 Let $2\le p, q\le\8$ and $p\not= q$. Let
$E\subset L_p(M)$ and $F\subset L_q(M)$ be closed subspaces. Let
$u: E\otimes F\to\comp$ be a linear functional and $c>0$ a
constant. The following assertions are equivalent
 \begin{enumerate}[{\rm i)}]
 \item  $u$ defines  a continuous functional on
$E\otimes_{h_{p,q}}F$ of norm $\le c$.
 \item  For all finite sequences $(a_k)\subset E$ and $(b_k)\subset F$
 \be
 \big|\sum_k u(a_k\otimes b_k)\big|\le c\,
 \big\|\big(\sum_ka_ka_k^*\big)^{1/2}\big\|_p\,
 \big\|\big(\sum_kb_k^*b_k\big)^{1/2}\big\|_q\;.
 \ee
 \item  There are positive unit functionals
$f\in \big(L_{p/2}(M)\big)^*$ and $g\in \big(L_{q/2}(M)\big)^*$
such that
 \be
 |u(a\otimes b)|\le c\, \big(f(aa^*)\, g(b^*b)\big)^{1/2}\ ,
 \quad a\in E,\ b\in F.
 \ee
 \item  There are positive unit functionals
$f\in \big(L_{p/2}(M)\big)^*$ and $g\in \big(L_{q/2}(M)\big)^*$
such that the associated linear map $\widetilde u: E\to F^*$
admits the following factorization
 \beq\label{dual of hpq4}
 \xymatrix{E \ar[d]_{i_E} \ar[r]^{\wt u} & F^*
 \\ H_p^r \ar[r]^{\wh u} & {\bar K}_q^r \ar[u]_{i_F^*}\;,}
 \eeq
where $i_E$ and $i_F$ are the natural inclusions associated to $f$
and $g$, respectively, given by Remark \ref{iE}, and where $\wh u$
is a \underbar{bounded} map.
\end{enumerate}
\end{rk}

In general, the map $\wh u$ in iv) above cannot be chosen to be
c.b.. Indeed, let $E=R_p$ and $F=C_q$ with $2\le p, q\le\8$. Then
both $R_p[R_p]$ and $C_q[C_q]$ are isometrically identified with
$S_2$. It follows that  $u: R_p\otimes C_q\to\comp$ satisfies ii)
above iff $\wt u: R_p \to C_{q'}$ is bounded. However, by using
the fact that $C_{q'}\cong R_q$, it is not hard to prove  that
$B(R_p, C_{q'})\not= CB(R_p, C_{q'})$ for $p\not= q$.

\begin{rk}\label{cs}
Let $E, F$ be two operator spaces and $u: E\times F\to\comp$  a
bilinear form. Let $p, q\ge 2$. Imitating the notion of completely
bounded bilinear forms in Christensen - Sinclair's sense (cf.
\cite{cs-rep}, \cite{cs-survey}), we say that $u$ is $(p,q)$-{\em
multiplicatively bounded} if there is a constant $c$ such that for
all $n\ge 1$ and all $a=(a_{i j})\in S_p^n[E], \ b=(b_{i j})\in
S_q^n[F]$
 $$\Big\|\Big(\sum_ku(a_{i k}, b_{k j})\Big)_{1\le i, j\le n}
 \Big\|_{S^n_{r}}\le c\,\|a\|_{S_p^n[E]}\,\|b\|_{S_q^n[F]}\; ,$$
where $1/r=1/p+1/q$. Let $\|u\|_{(p,q)-mb}$ denote the smallest of
such constants $c$. If $p=q$, $(p,q)$-multiplicatively bounded
forms are simply called $p$-{\em multiplicatively bounded} forms.
Now assume $E\subset L_p(M)$ and $F\subset L_q(M)$. Then $u$ is
$(p,q)$-multiplicatively bounded iff one of the assertions in
Remark \ref{dual of hpq} holds.
\end{rk}

\pf If $u$ is $(p,q)$-multiplicatively bounded, considering only
row and column matrices in the definition above, we see that the
assertion ii) of Remark \ref{dual of hpq} is verified. Conversely,
assume that iv) of Remark \ref{dual of hpq} holds. Let $\widetilde
u$ have the factorization (\ref{dual of hpq4}). Then for any $a\in
E$ and $b\in F$, $u(a, b)$ can be written as a product of three
operators:
 $$u(a, b)=i_E(a)\circ \hat u^*\circ i_F(b).$$
Recall that $i_F(b)\in K_q^c=S_q(\comp, K)$  and $i_E(a)\in
 H_p^r=S_p(\overline{H}, \comp)$.  Therefore, for any
$a=(a_{i j})\in S_p^n[E], \ b=(b_{i j})\in S_p^n[F]$
 $$\Big(\sum_ku(a_{i k}, b_{k j})\Big)_{i j}=\big[I_{S^n_p}\ot i_E(a)\big]
 \circ \big[I_{\el_2^n}\ot \hat u^*\big]
 \circ \big[I_{S^n_q}\ot i_F(b)\big].$$
Here $I_{S^n_q}\ot i_F(b)\in S_q^n[K_q^c]=S_q(\ell_2^n,
\ell_2^n(K))$, $I_{S^n_p}\ot i_E(a)\in S_p^n[H_p^r]=S_p(
\ell_2^n(\overline{H}), \ell_2^n)$ and $I_{\el_2^n}\ot \hat u^*\in
B(\ell_2^n(K), \ell_2^n(\overline{H}))$. Thus by the H\"older
inequality
 \be
 \Big\|\Big(\sum_ku(a_{i k}, b_{k j})\Big)_{i  j}
 \Big\|_{S^n_{r}}
 &=&\big\|\big[I_{S^n_p}\ot i_E(a)\big]
 \circ \big[I_{\el_2^n}\ot \hat u^*\big]
 \circ \big[I_{S^n_q}\ot i_F(b)\big]\big\|_{S^n_{r}}\\
 &\le&
 \big\|I_{S^n_p}\ot i_E(a)\big\|_{S_p^n[H_p^r]}\,
 \big\|I_{\el_2^n}\ot \hat u^*\big\|\,
 \big\|I_{S^n_q}\ot i_F(b)\big\|_{S_q^n[K_q^c]}\\
 &\le& \|\hat u\|\,\|a\|_{S_p^n[E]}\|b\|_{S_q^n[F]}\;.
 \ee
Therefore, $u$ is $(p, q)$-multiplicatively bounded, and so we
have proved the announced result. \cqd

\medskip

In particular, in the situation of Theorem \ref{dual of hp}, i.e.
when $p=q$ in Remark \ref{cs}, $u: E\times F\to \comp$ is
$p$-multiplicatively bounded iff the associated linear map $\wt u:
E\to F^*$ belongs to $\G_{R_p}(E, F^*)$. Consequently,
$p$-multiplicatively bounded forms are j.c.b.. Conversely, Theorem
\ref{G} implies that any j.c.b. form $u: E\times F\to \comp$
(still with $E, F\subset L_p(M)$ and $p\ge 2$) is the sum of a
$p$-multiplicatively bounded form and  the adjoint of a
$p$-multiplicatively bounded form. However, if $p\not=q$,  $(p,
q)$-multiplicatively bounded forms are in general not j.c.b..


\section{Noncommutative Khintchine inequalities}\label{khin}

In this section we give the main ingredient of the proofs of
Theorems \ref{G} and \ref{LG}. This is the noncommutative
Khintchine type inequalities for generalized circular systems. In
the sequel, $\H$ will be a fixed infinite dimensional separable
Hilbert space with an orthonormal basis $\{e_{\pm k}\}_{k\ge 1}$.
$\F(\H)$ stands for the associated free Fock space:
 $$\F(\H)=\bigoplus_{n=0}^\8 \H^{\otimes n}\ ,$$
where $\H^{\otimes 0}={\comp}\O$ with $\O$ a distinguished unit
vector. Let $\ell(e)$ (resp. $\ell^*(e)$) denote the left creation
(resp. annihilation) operator associated with a vector $e\in \H$.
Recall that $\ell^*(e)=(\ell(e))^*$.

Let $0\le \th\le 1$ and $\{\l_k\}_{k\ge 1}$ be a sequence of
positive numbers. Let
 \beq\label{circular}
 s_k=\el(e_{k}) + \l_k^{-1}\el^*(e_{-k})\quad\mbox{and}\quad
 g_k=\l_k^{\th}\,s_k,\quad k\ge 1.
 \eeq
The $s_k$ are generalized circular variables studied by
Shlyakhtenko \cite{shlya-quasifree}. We will also call
$(g_k)_{k\ge 1}$ a generalized circular system (with parameters
$\th$ and $(\l_k)$). Let $\G$ be the von Neumann algebra on
$\F(\H)$ generated by the $s_k$ (or equivalently by $g_k$). Let
$\rho$ be the vector state on $\G$ determined by the vacuum $\O$.
By \cite{shlya-quasifree}, $\rho$ is faithful on $\G$. Thus the
Haagerup $L_p$-space $L_p(\G)$ can be constructed from $\rho$. Let
$D$ denote the density of $\rho$ in $L_1(\G)$. Recall that $\rho$
can be recovered from $D$ as follows:
 $$\rho(x)=\tr(D x),\quad x\in \G.$$
Also recall that the modular group $\s_t^{\rho}$ is given by
 $$\s_t^{\rho}(x)=D^{it}\,xD^{-it},\quad
 x\in \G,\ t\in\real.$$
The $s_k$'s are eigenvectors of the modular group $\s_t^{\rho}$.
More precisely, we have the following formulas from
\cite{shlya-quasifree}
 \beq\label{eigenvector}
 \s_t(s_k)=\l_k^{-i2t}s_k\,,\quad \s_t(g_k)=\l_k^{-i2t}g_k\,,
 \quad k\ge 1,\ t\in
 \real
 \eeq
(see \cite[pp.342-343]{shlya-quasifree}; note that the minor
difference on parameters $\l$ between our definition of $s_k$
above and that of $y$ in \cite{shlya-quasifree}). The sequence
$\{g_k\}$ satisfies the following orthogonality with respect to
the state $\rho$: For any $0\le\eta\le1$
 \beq\label{circular orthogonal}
 \tr\big(g^*_j D^{\eta}g_kD^{1-\eta}\big)
 =\delta_{j,k}\;\l_k^{2(\th-\eta)}\ .
 \eeq
Indeed, by (\ref{eigenvector}), the left hand side of
(\ref{circular orthogonal}) is equal to
 $$\tr\big(g^*_j \s_{-i\eta}(g_k) D\big)
 =\tr\big(g^*_j \l_k^{-2\eta}g_k D\big)
 =\l_k^{-2\eta}\langle g_j\O, \ g_k\O\rangle
 =\delta_{j,k}\;\l_k^{2(\th-\eta)}\ .$$

\medskip

The following is the noncommutative  Khintchine type inequalities
for generalized circular systems. The case $p=\8$ was already
obtained in \cite{pisshlyak}.

\begin{thm}\label{khintchine}
Let  $1\le p\le\8$ and $\th=1/p$. Let $\{\l_k\}_k$ be a positive
sequence. Set
  \beq\label{circular p}
  g_{k,p}=D^{\frac{\th}{p}}\,g_k\,D^{\frac{1-\th}{p}}\;,
  \eeq
where $\{g_k\}_k$ is defined by $(\ref{circular})$. Let $M$ be a
von Neumann algebra  and $(x_n)$  a finite sequence in $L_p(M)$.
 \begin{enumerate}[{\rm i)}]
 \item  If $p\ge 2$,
 \begin{eqnarray}\label{k}
  \begin{array}{ccl}
 &&\begin{displaystyle}
 \max\Big\{
 \big\|\big(\sum_k\l_k^{2\th(1-\frac{2}{p})}\, x_k^*x_k
 \big)^{\frac{1}{2}}\big\|_p\;,\
 \big\|\big(\sum_k\l_k^{-2(1-\th)(1-\frac{2}{p})}\, x_kx_k^*
 \big)^{\frac{1}{2}}\big\|_p \Big\}
 \end{displaystyle}\\
 &&\hskip 3.5cm \begin{displaystyle}\le
 \big\|\sum_kx_k\ot g_{k, p}\big\|_{p}\le
 \end{displaystyle}\\
 && \begin{displaystyle}B_p\max\Big\{
 \big\|\big(\sum_k\l_k^{2\th(1-\frac{2}{p})}\, x_k^*x_k
 \big)^{\frac{1}{2}}\big\|_p\;,\
 \big\|\big(\sum_k\l_k^{-2(1-\th)(1-\frac{2}{p})}\, x_kx_k^*
 \big)^{\frac{1}{2}}\big\|_p \Big\}.\end{displaystyle}
  \end{array}
 \end{eqnarray}
 \item  If $p<2$,
 \begin{eqnarray}\label{kbis}
  \begin{array}{ccl}
 && \begin{displaystyle}A_p^{-1} \inf\Big\{
 \big\|\big(\sum_k\l_k^{2\th(1-\frac{2}{p})}\, a_k^*a_k
 \big)^{\frac{1}{2}}\big\|_p +
 \big\|\big(\sum_k\l_k^{-2(1-\th)(1-\frac{2}{p})}\, b_kb_k^*
 \big)^{\frac{1}{2}}\big\|_p\Big\}\end{displaystyle}\\
 &&\hskip 3.5cm \begin{displaystyle}\le
 \big\|\sum_kx_k\ot g_{k, p}\big\|_{p}\le\end{displaystyle}\\
 &&\begin{displaystyle}
 \inf\Big\{\big\|\big(\sum_k\l_k^{2\th(1-\frac{2}{p})}\, a_k^*a_k
 \big)^{\frac{1}{2}}\big\|_p +
 \big\|\big(\sum_k\l_k^{-2(1-\th)(1-\frac{2}{p})}\, b_kb_k^*
 \big)^{\frac{1}{2}}\big\|_p\Big\}.\end{displaystyle}
 \end{array}
 \end{eqnarray}
where the infimum runs over all decompositions $x_k=a_k+b_k$ in
$L_p(M)$. The two positive constants $A_p$ and $B_p$ depend only
on $p$ and can be controlled by a universal constant.
  \item  Let $G_{p}$ be the closed subspace of
$L_p(\G)$ generated  by $\{g_{k, p}\}_{k\ge1}$. Then there is a
completely bounded projection $P_{p}: L_p(\G)\to G_{p}$ such that
 $$\|P_{p}\|_{cb}\le 2^{|1-\frac{2}{p}|}\ .$$
 \end{enumerate}
\end{thm}

\n{\bf Remarks.} i) The  Khintchine inequalities above  also play
a crucial role in \cite{xu-emb} on the embedding of Pisier's OH
spaces, and more generally, the $q$-column spaces $C_q$  into
noncommutative $L_p$-spaces ($1\le p<q\le 2$).

 ii) (\ref{k}) is a particular case of a more general
inequality for free random series in \cite{jx-free}.

\medskip

{\em Proof of Theorem \ref{khintchine}} :  i) The state $\rho$ on
$\G$ extends to a contractive functional on $L_p(\G)$ for all
$1\le p\le\8$. More generally, let us consider the normal faithful
conditional expectation $\Phi\buildrel{\rm def}\over=I_M\ot\rho:
M\bar\ot\G\to M$. By \cite[Lemma 2.2]{jx-burk}, it extends to a
contractive projection from $L_p(M\bar\ot\G)$ onto $L_p(M)$ for
all $p\ge1$, still denoted by $\Phi$ in the following.

By \cite{shlya-quasifree}, the $g_k$'s are free in $(\G, \rho)$.
Thus by \cite{jx-free}, given a finite sequence $(x_k)\subset
L_p(M)$ $(2\le p\le\8)$ we have
 \be
 S\le \big\|\sum_k x_k\ot g_{k, p}\big\|_p\le
 B_p \,S,
 \ee
where
 \be
 S&=&\max\Big\{
 2^{1-\frac{2}{p}}\big(\sum_k\|x_k\ot g_{k, p}
 \|_p^p\big)^{\frac{1}{p}},\;
 \big\|\big(\sum_k \Phi(x_k^*x_k\ot g_{k,p}^*\,g_{k, p})
 \big)^{\frac{1}{2}}\big\|_p,\\
 && \hskip 5.3cm
 \big\|\big(\sum_k\Phi(x_kx_k^*\ot g_{k, p}\,g_{k, p}^*)
 \big)^{\frac{1}{2}}\big\|_p\Big\}.
 \ee
By (\ref{eigenvector}),
 $$g_{k, p}=\s_{-\frac{i\th}{p}}(g_k)\, D^{\frac{1}{p}}
 =\l_k^{-\frac{2\th}{p}}\, g_k\, D^{\frac{1}{p}}\;.$$
Thus,
 $$\Phi(x_k^*x_k\ot g_{k, p}^*\,g_{k, p})=
 x_k^*x_k\ot \big[\l_k^{-\frac{4\th}{p}}\,\rho(g_k^*\,g_k)\,
 D^{\frac{2}{p}}\big]
 =\l_k^{2\th(1-\frac{2}{p})}\,x_k^*x_k\ot
 D^{\frac{2}{p}}\;.$$
Therefore,
 $$\big\|\big(\sum_k \Phi(x_k^*x_k\ot g_{k, p}^*\,g_{k, p})
 \big)^{\frac{1}{2}}\big\|_p
 =\big\|\big(\sum_k\l_k^{2\th(1-\frac{2}{p})}\, x_k^*x_k
 \big)^{\frac{1}{2}}\big\|_p\;.$$
Similarly,
 $$\big\|\big(\sum_k \Phi(x_kx_k^*\ot g_{k, p}\,g_{k, p}^*)
 \big)^{\frac{1}{2}}\big\|_p
 =\big\|\big(\sum_k\l_k^{-2(1-\th)(1-\frac{2}{p})}\, x_kx_k^*
 \big)^{\frac{1}{2}}\big\|_p\;.$$
Combining the previous inequalities, we get the lower estimate in
(\ref{k}).

To prove the upper estimate we only need to show that the term
$(\sum_k\|x_k\ot g_{k, p}\|_p^p)^{1/p}$ is controlled by the two
others. To this end we first observe that
 $$\|g_{k, p}\|_{p}\le \|g_{k}\|_{\8}^{1-\frac{2}{p}}\,
 \|g_{k,2}\|_{2}^{\frac{2}{p}}
 \le (\l_k^\th + \l_k^{-1+\th})^{1-\frac{2}{p}}\;;$$
whence
 $$\|x_k\ot g_{k, p}\|_{p}\le
 (\l_k^\th + \l_k^{-1+\th})^{1-\frac{2}{p}}\,\|x_k\|_p\ .$$
Thus
 $$\big(\sum_k\|x_k\ot g_{k, p}\|^p_{p}\big)^{\frac{1}{p}}\le
 \big(\sum_k\l_k^{\th(1-\frac{2}{p})}\,\|x_k\|_{p}
 \big)^p\big)^{\frac{1}{p}}
 + \big(\sum_k\l_k^{-(1-\th)(1-\frac{2}{p})}\,
 \|x_k\|_{p}\big)^p\big)^{\frac{1}{p}}\ .$$
However, for any $(y_k)\subset L_p(M)$ we have
 $$\big(\sum_k \|y_k\|_{p}^p\big)^{\frac{1}{p}}\le\min\big\{
 \big\|\big(\sum_k y_k^*y_k\big)^{\frac{1}{2}}\big\|_{p}\ ,\
 \big\|\big(\sum_k y_ky_k^*\big)^{\frac{1}{2}}\big\|_{p}\big\}\ .$$
Therefore,
 $$\big(\sum_k\|x_k\ot g_{k, p}\|^p_{p}\big)^{\frac{1}{p}}\le
 \big\|\big(\sum_k\l_k^{2\th(1-\frac{2}{p})}\, x_k^*x_k
 \big)^{\frac{1}{2}}\big\|_p +
 \big\|\big(\sum_k\l_k^{-2(1-\th)(1-\frac{2}{p})}\,
 x_kx^*_k\big)^{\frac{1}{2}}\big\|_p\;;$$
whence the upper estimate in ( \ref{k}).

\medskip
{\rm ii)} The minoration here follows from the majoration in {\rm
i)} by a simple duality argument. Indeed, let $(y_k)\subset
L_{p'}(M)$ be a finite sequence such that
 $$\max\Big\{
 \big\|\big(\sum_k\l_k^{2\th(1-\frac{2}{p'})}\, y_k^*y_k
 \big)^{\frac{1}{2}}\big\|_{p'}\;,\
 \big\|\big(\sum_k\l_k^{-2(1-\th)(1-\frac{2}{p'})}\, y_ky_k^*
 \big)^{\frac{1}{2}}\big\|_{p'} \Big\}\le 1.$$
Let
 $$x=\sum_kx_k\ot g_{k, p}\quad\mbox{and}\quad
   y=\sum_k y_k\ot g_{k, p'}\;.$$
Then by (\ref{circular orthogonal})
 $$\sum_k\tr(y_k^*x_k)=\tr\ot\tr(y^*x)$$
and by {\rm i)}
 \begin{eqnarray*}
 \|y\|_{p'}\le B_{p'}
 \max\Big\{
 \big\|\big(\sum_k\l_k^{2\th(1-\frac{2}{p'})}\, y_k^*y_k
 \big)^{\frac{1}{2}}\big\|_{p'}\;,\
 \big\|\big(\sum_k\l_k^{-2(1-\th)(1-\frac{2}{p'})}\, y_ky_k^*
 \big)^{\frac{1}{2}}\big\|_{p'} \Big\}\le B_{p'}.
 \end{eqnarray*}
Hence
 $$\big|\sum_k\tr(y_k^*x_k)\big|
 \le \|x\|_p\|y\|_{p'}\le B_{p'}\|x\|_p\;.$$
Taking the supremum over all $(y_k)$ as above yields the lower
estimate in (\ref{kbis}) with $A_p=B_{p'}$. The majoration is a
consequence of the following elementary inequality (with $1\le
p\le 2$)
 $$\|x\|_p\le \big\|\big[\Phi(x^*x)\big]^{\frac{1}{2}}\big\|_p\;,
 \quad\forall\;x\in L_p(M\bar\ot\G).$$
By duality, this immediately follows from
 $$\|\Phi(y^*y)\|_{p'/2}\le \|y^*y\|_{p'/2}\;,\quad
 \forall\;y\in L_{p'}(M\bar\ot\G).$$

\medskip

iii) The following type of arguments is rather standard today (cf.
\cite{hp-duke} in the case of free groups).  Recall that $\O$ is a
separating vector for $\G$ (cf. \cite{shlya-quasifree}). Thus any
operator $a\in \G$ is uniquely determined by $a\O$. Set
$\G_2=\{a\O\ :\ a\in\G\}$. Then $\G_2$ is a vector subspace of
$\F(\H)$, which is isometric to $L_2(\G)$. For any $\xi\in\G_2$ we
denote by $W(\xi)$ the unique operator in $\G$ such that
$W(\xi)\O=\xi$. It is well known (and easy to check) that all
tensors from $\H^{\otimes n}$ belong to $\G_2$ ($n\in\nat$). Now
we use multi-index notation. Recall that $\{e_{\pm k}\}_{k\ge1}$
is an orthonormal basis of $\H$. For any $i_1,...,
i_n\in\ent\setminus\{0\}$ we put $\underline i=(i_1,..., i_n)$ and
$e_{\underline i} = e_{i_1}\otimes ... \otimes e_{i_n}$. If
$\underline i=\emptyset$, we set $e_{\underline i}=\O$. Then
$\{e_{\underline i}\}_{\underline i}$ is an orthonormal basis of
$\F(\H)$. By the discussion above, every $e_{\underline i}$
belongs to $\G_2$ and $W(e_{\underline i})$ is the unique operator
in $\G$ such that $W(e_{\underline i})\O=e_{\underline i}$. By the
definition of $g_k$ in (\ref{circular}), we have
 $$W(e_{k})=\l_k^{-\th}\,g_k\quad \mbox{and}
 \quad W(e_{-k})=\l_k^{1-\th}\,g^*_k\ ,\quad k\ge 1.$$
Given $a\in\G$, developing $a\O$ in the orthonormal basis
$\{e_{\underline i}\}_{\underline i}$ we can (symbolically) write
 $$a=\sum_{\underline i}\ c_{\underline i}(a)\;W(e_{\underline i})\ ,$$
where $c_{\underline i}= \la e_{\underline i},\; a\O\ra$. Let
$\G_0$ be the $\ast$-subalgebra of $\G$ of all operators which
admit a finite development as above. Note that $\G_0$ is w*-dense
in $\G$. Consequently, $D^{1/2p}\,\G_0\,D^{1/2p}$ is dense in
$L_p(\G)$ for any $1\le p<\8$ (cf. \cite[Lemma 1.1]{jx-burk}; see
also \cite{ju-doob}).

Now for any $a\in \G_0$ we define
 $$P_\8(a)=\sum_{k\ge 1}\ c_{k}(a)\;W(e_{k})$$
and for $1\le p<\8$
 $$P_p : D^{\frac{1}{2p}}\,\G_0\,D^{\frac{1}{2p}}\longrightarrow
 D^{\frac{1}{2p}}\,\G_0\,D^{\frac{1}{2p}}\quad\mbox{by}\quad
 P_p\big(D^{\frac{1}{2p}}\,a\,D^{\frac{1}{2p}}\big)
 = D^{\frac{1}{2p}}\,P_\8(a)\,D^{\frac{1}{2p}}\ .$$
We are going to show that $P_p$ extends to a completely bounded
projection from $L_p(\G)$ onto $G_p$. In fact, what will be shown
is that $I\ot P_p$ extends to a bounded projection from
$L_p(M\bar\ot\G)$ onto $L_p(M)\bar\ot G_p$ for any von Neumann
algebra $M$, where $L_p(M)\bar\ot G_p$ is the closure of $L_p(M)
\ot G_p$ in $L_p(M\bar\ot\G)$ (relative to the w*-topology for
$p=\8$). This is clear for $p=2$; moreover, the extension of $I\ot
P_2$ is the orthogonal projection from $L_2(M\bar\ot\G)$ onto
$L_2(M)\bar\ot G_2$.

Then consider the case $p=\8$. Let $\{x_{\underline
i}\}_{\underline i}$ be a finite family in $M$. Let
 $$x=\sum_{\underline i}\ x_{\underline i}\otimes
 W(e_{\underline i})\in M\ot\G. $$
Then
 $$I\ot P_\8(x)=\sum_{k\ge 1}\ x_{k}\otimes W(e_{k})
 =\sum_{k}\ \l_k^{-\th}\;x_{k}\ot g_k .$$
Therefore, by (\ref{k})
 $$\|I\ot P_\8(x)\|_\8\le 2
 \max\big\{
\big\|\big(\sum_k\ x_{k}^*\,x_{k}\big)^{1/2}\big\|_\8\ ,\
 \big\|\big(\sum_k\l_k^{-2}\, x_{k}\,x^*_{k}
 \big)^{1/2}\big\|_\8 \big\}\ .$$
Let $\xi\in H$ be a unit vector ($H$ being the Hilbert space at
which $M$ acts). Then by the orthonormality of $\{e_{\underline
i}\}_{\underline i}$
 \be
 \|x\|_\8^2
 &\ge& \langle x(\xi\otimes\O),\ x(\xi\otimes\O)\rangle
 =\sum_{\underline i}\|x_{\underline i}(\xi)\|^2\\
 &=&\langle\xi,\ \sum_{\underline i}\ x_{\underline i}^*\,
 x_{\underline i}(\xi)\rangle\ge
 \langle\xi,\ \sum_k\ x_{k}^*\,x_{k}(\xi)\rangle ;
 \ee
whence
 $$\|\sum_k\ x_k^*x_k\|_\8\le \|x\|_\8^2.$$
On the other hand,
 \be
 \langle x^*(\xi\otimes\O),\ x^*(\xi\otimes\O)\rangle
 &=&\sum_{{\underline i},\, \underline j}
 \la x_{\underline i}^*(\xi),\ x_{\underline j}^*(\xi)\ra\;
 \la W(e_{\underline i})^*(\O),\ W(e_{\underline j})^*(\O)\ra\\
 &=&\sum_{\underline i,\, \underline j}
 \la \xi,\ x_{\underline i}\;x_{\underline j}^*(\xi)\ra\;
 \rho\big[W(e_{\underline i})  W(e_{\underline j})^*\big]\\
 &=&\sum_{\underline i,\, \underline j}
 \la \xi,\ x_{\underline i}\;x_{\underline j}^*(\xi)\ra\;
 \rho\big[W(e_{\underline j})^*\s_{-i}(W(e_{\underline
 i}))\big].
 \ee
However, by \cite{shlya-quasifree} one can easily show
 $$\s_{-i}\big(W(e_{\underline i})\big)
 =\b_{\underline i}\;W(e_{\underline i}),$$
where $\b_{\underline i}$ is a finite product involving
$\l_k^{\pm\th}$ and $\l_k^{\pm(1-\th)}$. Then it follows that
 $$\rho\big[W(e_{\underline j})^*\,\s_{-i}(
 W(e_{\underline i}))\big]
 =\la W(e_{\underline j})\O, \ \s_{-i}
 \big(W(e_{\underline i})\big)\O\ra
 =\d_{\underline i,\;\underline j}\;\b_{\underline i}^2.$$
Hence
 \be
 \langle x^*(\xi\otimes\O),\ x^*(\xi\otimes\O)\rangle
 &\ge& \sum_{k\ge1}\la \xi,\ x_{k}x_{k}^*(\xi)\ra\;
 \rho\big[W(e_{k})^*\,\s_{-i}(W(e_k))\big]\\
 &=&\sum_{k\ge1}\la \xi,\ x_{k}x_{k}^*(\xi)\ra\;\l_k^{-2}
 =\la \xi,\ \sum_{k\ge1}\l_k^{-2}\,x_{k}x_{k}^*(\xi)\ra .
 \ee
Therefore
 $$\big\|\sum_{k\ge1}\l_k^{-2}\,x_{k}x_{k}^*\big\|_\8
 \le \|x\|_\8^2.$$
Combining the preceding inequalities, we obtain
 $$\|I\ot P_\8(x)\|_\8\le 2\|x\|_\8.$$
This is the key inequality of this part of the proof. From this we
get the extension property for $p=1$ by virtue of the following
easily checked duality equality: for any $a, b\in\G_0$ and $1\le
p\le\8$
 $$\la P_p\big(D^{\frac{1}{2p}}\,b\,D^{\frac{1}{2p}}\big),\
 D^{\frac{1}{2p'}}\,a\,D^{\frac{1}{2p'}}\ra
 =\la D^{\frac{1}{2p}}\,b\,D^{\frac{1}{2p}},\
 P_{p'}\big(D^{\frac{1}{2p'}}\,a\,D^{\frac{1}{2p'}}\big)\ra.$$
Indeed, by this equality combined with the preceding boundedness
of $P_\8$ on $\G_0$, we deduce that $I\ot P_1$ is bounded on
$L_1(M)\ot\big[D^{\frac{1}{2}}\,\G_0\, D^{\frac{1}{2}}\big]$ with
respect to the $L_1$-norm and is of norm $\le 2$. Thus by the
density of $D^{\frac{1}{2}}\,\G_0\, D^{\frac{1}{2}}$ in $L_1(\G)$,
we deduce that $I\ot P_1$ extends to a bounded map on
$L_1(M\bar\ot\G)$. By duality once more, we see that the adjoint
of this extension of $I\ot P_1$ yields the desired (normal)
extension of $I\ot P_\8$ on $M\bar\ot\G$. The remaining case for
$1<p<2$ or $2<p<\8$ is proved by Kosaki's interpolation theorem
\cite{ko-int}. Therefore the proof of Theorem \ref{khintchine} is
complete.\cqd


\section{Vector-valued noncommutative $L_p$}\label{vvLp}

This section contains the second main tool of the proofs of
Theorems \ref{G} and \ref{LG}, i.e. the vector-valued
noncommutative $L_p$-spaces for QWEP von Neumann algebras.  The
theory of vector-valued noncommutative $L_p$-spaces was first
developed by Pisier \cite{pis-asterisque} for injective semifinite
von Neumann algebras. Very recently, Junge \cite{ju-fubini},
\cite{ju-qwep} partly extended this theory to QWEP von Neumann
algebras. Junge's idea is to represent QWEP von Neumann algebras
as images of normal conditional expectations on ultraproducts of
injective von Neumann algebras, and then apply Pisier's theory.
His approach relies heavily upon the theory of ultraproducts of
noncommutative $L_p$-spaces developed recently by Raynaud
\cite{ray-ultra}.

We first recall some known results on the ultraproducts  of von
Neumann algebras and noncommutative $L_p$-spaces. Let $\U$ be a
free ultrafilter on some index set $I$. If $X$ is an operator
space, we use the notation $X^\U$ to denote the ultrapower of $X$.
$X^\U$ is equipped with its natural operator space structure as
introduced by Pisier (cf. \cite{pis-intro}). Now let $M$ be a von
Neumann algebra. Groh proved that the ultrapower $(M_*)^\U$ of the
predual $M_*$ is again a predual of  von Neumann algebra (cf.
\cite{groh}; see also \cite{ray-ultra}). In fact, assuming $M$
acts standardly on some Hilbert space $H$, the ultrapower $M^\U$
is a C*-algebra, which can be naturally represented on the Hilbert
space ultrapower $H^\U$. Then the von Neumann algebra
$\big((M_*)^\U\big)^*$ is the w*-closure of $M^\U$ in $B(H^\U)$.
In the sequel we will denote $\big((M_*)^\U\big)^*$ by $M_\U$. On
the other hand,  Raynaud developed the theory of ultraproducts of
noncommutative $L_p$-spaces. In particular, he proved that the
ultrapower $\big(L_p(M)\big)^\U$ can be identified with
$L_p(M_\U)$; moreover, this identification is natural in the sense
that it preserves all algebraic operations such as product,
involution, positivity ...

\medskip

Recall that a C*-algebra is called WEP (for weak expectation
property) in Lance's sense if the natural inclusion
$A\hookrightarrow A^{**}$ can be factorized completely
contractively through some $B(H)$. $A$ is called QWEP if $A$ is a
quotient of a WEP C*-algebra (cf. \cite{kir-nonsemi}).

The following characterization of QWEP due to  Junge
\cite{ju-fubini} will play an important role later.

\begin{prop}\label{qwep}
 A von Neumann algebra $M$ is QWEP iff there are a Hilbert
space $H$ and a free ultrafilter $\U$ on some index set $I$ such
that $M$ is the image of a normal conditional expectation on
$B(H)_\U$.
\end{prop}

\medskip
Let $M$ and $B(H)_\U$ be as in the proposition above. Let $\Phi:
B(H)_\U\to M$ be the corresponding normal conditional expectation
($\Phi$ is, in general, not faithful). In this case, $L_p(M)$ can
be naturally identified as a subspace of $L_p(B(H)_\U)$. It is
also know that $\Phi$ defines a contractive projection $\Phi_p$
from $L_p(B(H)_\U)$ onto $L_p(M)$ for any $1\le p\le\8$ (cf.
\cite[Proposition~2.3] {jx-burk}). Let $N$ be any von Neumann
algebra. Then $I_{L_p(N)}\otimes \Phi_p$ is also a contractive
projection from $L_p(N\overline{\otimes} B(H)_\U)$ onto
$L_p(N\overline{\otimes} M)$. Indeed, it is obvious that
$I_{N}\otimes \Phi$ is a normal conditional expectation from
$N\overline{\otimes} B(H)_\U$ onto $N\overline{\otimes} M$.
Applying the previous result to $I_{N}\otimes \Phi$, we get the
announced one. In particular, $\Phi_p$ is completely contractive.

\medskip

Now we introduce the vector-valued noncommutative $L_p$-spaces for
QWEP von Neumann algebras. Let $M$ and $H$ be as in Proposition
\ref{qwep}. Given an operator space $X$, the $X$-valued Schatten
class $S_p[H; X]$ is defined  in section 1. Now let $E$ be a
finite dimensional operator space. Following Junge, we define
$L_p[M; E]$ ($1\le p<\8$) simply  as the ultrapower $\big(S_p[H;
E]\big)^\U$. Then for any operator space $X$, $L_p[M; X]$ is
defined as the closure in $\big(S_p[H; X]\big)^\U$ of $L_p[M; E]$
when $E$ runs over all finite dimensional spaces of $X$.

\medskip

We will need the following from \cite{ju-qwep}.

\begin{prop}\label{cb extension}
 Let $M$ be QWEP and $1\le p<\8$.
\begin{enumerate}[{\rm i)}]
 \item  Let $u: E\to F$ be a c.b. map between operator spaces.
Then  $I_{L_p(M)}\otimes u$ extends to a c.b. map from $L_p[M; E]$
into $L_p[M; F]$ and $\|I_{L_p(M)}\otimes u\|_{cb}\le \|u\|_{cb}$.
Moreover, if $u$ is a complete isometry, then so is
$I_{L_p(M)}\otimes u$.
 \item  Let $N$ be another von Neumann algebra. Then
  $$L_p[M; L_p(N)]=L_p(M\overline{\otimes} N)\quad
  \mbox{completely isometrically}.$$
\end{enumerate}
\end{prop}

Note that part i) above easily follows from the definition and the
corresponding results on the vector-valued Schatten classes in
\cite{pis-asterisque}. Part ii) is more substantial. It is a
consequence of Junge's noncommutative Fubini theorem (cf.
\cite{ju-fubini}). We refer to \cite{ju-fubini} and \cite{ju-qwep}
for more details.

\begin{prop}\label{p to q}
 Let $E$ be an operator space, $1\le p<q\le\8$ and
$1/r=1/p-1/q$.
\begin{enumerate}[{\rm i)}]
 \item  If $q<\8$, then for any
$x\in L_q(B(H)_\U)\otimes E$
 \beq\label{p to q1}
 \|x\|_{L_q[B(H)_\U;E]}=\sup\big\{
 \|axb\|_{L_p[B(H)_\U;E]} \big\},
 \eeq
where the supremum runs over all $a$ and $b$ in the unit ball of
$L_{2r}(B(H)_\U)$.
 \item  Assume $q=\8$ $($so $r=p)$. If $E$ is exact, then for any
 $x\in  B(H)^\U\otimes E$
 \beq\label{p to q2}
 \|x\|_{B(H)^\U\otimes_{\rm min}E}
 \le\sup\big\{\|axb\|_{L_p[B(H)_\U;E]}\big\}
 \le\l\, \|x\|_{B(H)^\U\otimes_{\rm min}E}\;,
 \eeq
where the supremum runs over all $a$ and $b$ in the unit ball of
$L_{2p}(B(H)_\U)$, and where $\l=ex(E)$ is the exactness constant
of $E$. Conversely, if (\ref{p to q2}) holds for some infinite
dimensional Hilbert space $H$ and some constant $\l$, then $E$ is
exact and $ex(E)\le\l$.
\end{enumerate}
\end{prop}

\pf i) Let $x\in L_q(B(H)_\U)\otimes E$. Passing to a finite
dimensional subspace if necessary, we may assume $E$ finite
dimensional. Then by definition
 $$L_q[B(H)_\U;E]=\big(S_q[H; E]\big)^\U\ .$$
By (1.5), for any $x_i\in S_q[H; E]$ we have
 $$\|x_i\|_{S_q[H; E]}=\sup\big\{
 \|a_ix_ib_i\|_{S_p[H; E]}  :  a_i,  b_i\in S_{2r}(H),
 \|a_i\|_{2r}\le 1,\;\|b_i\|_{2r}\le 1\big\}.$$
On the other hand, by [Ra] for any $t<\8$
 $$L_t(B(H)_\U)=\big(S_t(H)\big)^\U\ .$$
Combining these, we easily deduce (\ref{p to q1}).

ii) Again, we can assume $E$ finite dimensional. By (1.4), for any
$x_i\in B(H)\otimes E$
 $$\|x_i\|_{B(H)\otimes_{\rm min}E}=
 \sup\big\{\|a_ix_ib_i\|_{S_p[H; E]}  :  a_i,
 b_i\in S_{2p}(H), \|a_i\|_{2p}\le 1, \|b_i\|_{2p}\le 1\big\}.$$
Thus as before, we get
 $$\|x\|_{\big(B(H)\otimes_{\rm min}E\big)^\U}=
 \sup\big\{\|axb\|_{L_p[B(H)_\U;E]} : a,
 b\in L_{2p}(B(H)_\U), \|a\|_{2p}\le 1,  \|b\|_{2p}\le 1\big\}\ .$$
Therefore (\ref{p to q2}) can be rewritten as
 $$\|x\|_{B(H)^\U\otimes_{\rm min}E}\le
 \|x\|_{\big(B(H)\otimes_{\rm min}E\big)^\U}
 \le \l\, \|x\|_{B(H)^\U\otimes_{\rm min}E}\ .$$
It is known that  the first inequality above is always true, while
the validity of the second is equivalent to the exactness of $E$;
moreover, the least constant $\l$ is then equal to the exactness
constant of $E$. We omit the details and refer to \cite[Chapiter
17] {pis-intro}. \cqd

\begin{cor}\label{multiplication map}
 Let $p, q, r$ be as in Proposition \ref{p to q}. Given $a, b\in
L_{2r}(B(H)_\U)$ we define  $M_{a, b}(x)=axb$.
\begin{enumerate}[\rm i)]
 \item  If $q<\8$, $M_{a, b}$ defines a c.b. map from
$L_q[B(H)_\U;E]$ into $L_p[B(H)_\U;E]$ and $\|M_{a, b}\|_{cb}\le
\|a\|_{2r}\|b\|_{2r}$.
 \item  If $q=\8$ and $E$ is exact, $M_{a, b}$ defines a c.b. map from
$\big(B(H)\big)^\U\otimes_{\rm min}E$ into $L_p[B(H)_\U;E]$ and
$\|M_{a, b}\|_{cb}\le ex(E)\|a\|_{2p}\|b\|_{2p}$.
\end{enumerate}
\end{cor}

\pf We only prove i). The proof for ii) is similar. It is
immediate from Proposition \ref{p to q} that $M_{a, b}$ is bounded
and $\|M_{a, b}\|\le \|a\|_{2r}\|b\|_{2r}$. (Only this boundedness
will be needed later.) To prove the complete boundedness we use
Lemma \ref{cb via Sp}, so we have to show
 $$\big\|I_{S^n_q}\ot M_{a, b}: S^n_q\big[L_q[B(H)_\U;E]\big]\to
 S^n_q\big[L_p[B(H)_\U;E]\big]\big\|\le \|a\|_{2r}\|b\|_{2r}\;,
 \quad \forall\;n\in\nat.$$
However,
 $$S^n_q\big[L_q[B(H)_\U;E]\big]=L_q[B(\el^n_2(H))_\U;E].$$
Thus for any $\a, \b\in S^n_{2r}$, $M_{\a\ot a,\, \b\ot b}$ is
bounded from $L_q[B(\el^n_2(H))_\U;E]$ to
$L_p[B(\el^n_2(H))_\U;E]$ and of norm $\le \|\a\ot a\|_{2r}\|\b\ot
b\|_{2r}= \|\a\|_{2r}\|a\|_{2r}\|\b\|_{2r}\|b\|_{2r}$. Taking the
supremum over all $\a$ and $\b$ such that $\|\a\|_{2r}\le 1$ and
$\|\b\|_{2r}\le 1$, and using (\ref{Sq via Sp}), we deduce the
announced result. \cqd

\section{Proof  of Theorem \ref{G}}\label{pfG}

This section and the next are devoted to the proofs of Theorems
\ref{G} and \ref{LG}, respectively. The common key ingredient of
both proofs  is Theorem \ref{khintchine}. The patterns of our
proofs are similar to those of the corresponding results for
$L_\8$ in \cite{pisshlyak} and \cite{pis-Rp}.

\medskip

\n{\it Proof of Theorem \ref{G}.} i) $\Rightarrow$ ii). Let
$(\mu_k)_k$ be a sequence of positive numbers. Set
$\l_k=\mu_k^{p/(p-2)}$. Let $\{g_k\}_k$ be the generalized
circular system with parameters $(\l_k)_k$ and $\th=1/2$ defined
by (\ref{circular}). We will use the notations introduced in
section~3: $\G$ is the von Neumann algebra generated by the $g_k$,
$\rho$ the normal faithful state on $\G$ given by the vacuum  and
$D$ the density of $\rho$ in $L_1(\G)$. It was proved in
\cite{pisshlyak} that $\G$ is QWEP. Thus by Proposition
\ref{qwep}, $\G$ is the image of a normal conditional expectation
$\Phi$ on some $B(H)_\U$. Consequently, $L_p(\G)$ can be naturally
identified as a subspace of $L_p(B(H)_\U)$.

To prove ii) we can clearly assume $E$ and $F$ finite dimensional.
Then by definition, $L_p[\G; E]$ is a subspace of $L_p[B(H)_\U;
E]$. By Proposition \ref{cb extension},
 $$v\mathop{=}^{\rm def} I_{L_p(B(H)_\U)}\otimes \widetilde u\ :\
 L_p[B(H)_\U; E]\to L_p[B(H)_\U; F^*]\quad \mbox{is bounded}.$$
On the other hand, letting $\frac{1}{r}=\frac{1}{p'}-\frac{1}{
p}=1-\frac{2}{p}$ (so $r$ is the conjugate index of $p/2$) and by
Corollary \ref{multiplication map} (noting that
$\|D^{1/2r}\|_{2r}=1$),
 $$w\mathop{=}^{\rm def} M_{D^{1/2r}, D^{1/2r}}\ :\ L_p[B(H)_\U; F^*]\to
 L_{p'}[B(H)_\U; F^*]\quad \mbox{is contractive}.$$
Therefore,
 $$\big\|wv\; :\; L_p[B(H)_\U; E]\to L_{p'}[B(H)_\U; F^*]
 \big\|\le\|u\|_{jcb}\ .$$
Now let $(a_k)\subset E$  be a finite sequence. Then (recalling
that $g_{k, p}$ is defined by (\ref{circular p}) with
$\th=\frac{1}{2}$)
 $$wv\Big(\sum_k a_k\ot g_{k, p}\Big)
 =\sum_k \widetilde u(a_k)\ot g_{k, p'}\ .$$
Thus
 \be
 \Big\|\sum_k
  \widetilde u(a_k)\ot
 g_{k, p'}\Big\|_{L_{p'}[B(H)_\U; F^*]}
 \le \|u\|_{jcb}\,\Big\|\sum_k a_k\ot
 g_{k, p}\Big\|_{L_p[B(H)_\U; E]} \;.
 \ee
Let $F^\perp\subset L_{p'}(M)$ be the orthogonal complement of
$F$. Then
 $$L_{p'}[B(H)_\U; F^*]=\big(S_{p'}[H; F^*]\big)^\U
 =\frac{\big(S_{p'}[H; L_{p'}(M)]\big)^\U}
 {\big(S_{p'}[H; F^\perp]\big)^\U}\;.$$
On the other hand, by Proposition \ref{cb extension}, we have the
following isometric inclusions
 $$L_t[B(H)_\U; X]\subset L_t[B(H)_\U; L_t(M)]=
 L_t(B(H)_\U\;\overline{\otimes}\;M)\subset \big(S_{t}[H;
 L_{t}(M)]\big)^\U$$
for any $1\le t<\8$ and any subspace $X\subset L_t(M)$. It is
clear that
 $$\big(S_{p'}[H; F^\perp]\big)^\U\subset
 \big(L_p[B(H)_\U; F]\big)^\perp\ .$$

Given a von Neumann algebra $N$ we use the following duality
bracket between $L_p(N)$ and $L_{p'}(N)$ in the category of
operator spaces
 $$\la y, x\ra=\tr(\overline{y^*}\,x),\quad x\in L_p(N), \;
 y\in L_{p'}(N).$$
This duality is consistent with the operator space structure on
$L_p(N)$ (recalling that $L_1(N)$ is the predual of $N^{op}$ and
$N^{op}\cong\overline N$). With this duality, the dual of $L_p(N)$
is $\overline{L_{p'}(N)}$ completely isometrically.

 Note that by
(\ref{circular orthogonal}) we have
 $$\la g_{j, p'},\; g_{k, p}\ra
 =\tr(g_j^*D^{\frac{1}{2}}g_kD^{\frac{1}{2}})=\d_{j,k}\;.$$
Then for any finite sequence $(b_k)\subset F$ we deduce  that
 \be
 \big|\sum_k\;u(a_k, b_k)\big|
 &=&\big|\sum_k\;\widetilde u(a_k)(b_k)\big|\\
 &=&\Big|\big\langle\sum_k\widetilde u(a_k)\ot g_{k, p'}\;,\
 \sum_kb_k\ot g_{k, p}\big\rangle\Big|\\
 &\le& \Big\|\sum_kb_k\ot g_{k, p}
 \Big\|_{L_p[B(H)_\U;F]}\,\Big\|\sum_k\widetilde u(a_k)\ot g_{k, p'}
 \Big\|_{L_{p'}[B(H)_\U; F^*]}\;.
 \ee
Therefore, combining the preceding inequalities with Theorem
\ref{khintchine}, we get
 \be
\big|\sum_k\;u(a_k, b_k)\big|
 &\le& \|u\|_{jcb}\Big\|\sum_k a_k\ot g_{k, p}
 \Big\|_{L_p[B(H)_\U; E]}\,\Big\|\sum_kb_k\ot g_{k, p}
  \Big\|_{L_p[B(H)_\U;F]}\\
 &=& \|u\|_{jcb}\Big\|\sum_k a_k\ot g_{k, p}
 \Big\|_{L_p(\G\overline{\otimes}M)}\,\Big\|\sum_kb_k\ot g_{k, p}
 \Big\|_{L_p(\G\overline{\otimes} M)}\\
 &\le& B_p\|u\|_{jcb}\Big[\Big\|(\sum_k\l_k^{1-\frac{2}{p}}a_k^*
 a_k)^{\frac{1}{2}}\Big\|_p
 + \Big\|(\sum_k\l_k^{-(1-\frac{2}{p})}a_ka_k^*)^{\frac{1}{2}}\Big\|_p\Big]\\
 &&\hskip 1.4cm\Big[\Big\|(\sum_k\l_k^{1-\frac{2}{p}}b_k^*
 b_k)^{\frac{1}{2}}\Big\|_p
 + \Big\|(\sum_k\l_k^{-(1-\frac{2}{p})}b_kb_k^*)^{\frac{1}{2}}\Big\|_p\Big].
 \ee
This is (\ref{G2}) by the relation between $\l_k$ and $\mu_k$.

\medskip

 ii) $\Rightarrow$ iii). This is done by a standard Hahn-Banach separation
argument as in \cite{pisshlyak}. For completeness, we include the
main lines. Assume all the $\mu_k$ are equal, say,  to $s$. Then
by (\ref{G2})
 \be
 \Big|\sum_ku(a_k, b_k)\Big|
 &\le& 2K_2
 \Big[\Big\|\sum_k\,s\,a_k^*a_k\Big\|_{\frac{p}{2}}
 + \Big\|\sum_k\,s^{-1}a_ka_k^*\Big\|_{\frac{p}{2}}\Big]^{\frac{1}{2}}\\
 &&\hskip 0.7cm\Big[\Big\|\sum_k\,s\,b_k^*b_k\Big\|_{\frac{p}2}
 + \Big\|\sum_k\,s^{-1}b_kb_k^*\Big\|_{\frac{p}{2}}\Big]^{\frac{1}{2}}\\
 &\le& K_2
 \Big[\Big\|\sum_k\,s\,a_k^*a_k\Big\|_{\frac{p}{2}}
 + \Big\|\sum_k\,s^{-1}a_ka_k^*\Big\|_{\frac{p}{2}}\\
 &&\hskip 0.4cm+\Big\|\sum_k\,s\,b_k^*b_k\Big\|_{\frac{p}{2}}
 + \Big\|\sum_k\,s^{-1} b_kb_k^*\Big\|_{\frac{p}{2}}\Big].
 \ee
Then by a Hahn-Banach argument as in the proof of Proposition
\ref{to Rp}, we get positive operators $f_1, f_2, g_1, g_2\in
L_r(M)$ ($r$ being the conjugate index of $p/2$), all of them with
norms $\le 1$, such that for any $a\in E,\ b\in F$
 $$|u(a, b)|\le K_2\big[f_1(s^{-1}aa^*)+ f_2(sa^*a)+
 g_1(sb^*b) +g_2(s^{-1}bb^*)\big].$$
Replacing $a$ and $b$ respectively by $ta$ and $t^{-1}b$ in the
above inequality and then taking the infimum over all $t>0$, we
deduce that
 \be
 |u(a, b)|
 &\le& 2K_2\big[f_1(s^{-1}aa^*)+f_2(sa^*a)\big]^{1/2}
 \big[g_1(sb^*b) +g_2(s^{-1}bb^*)\big]^{1/2}\\
 &=&2K_2\big[f_1(aa^*)g_1(b^*b)
 +f_2(a^*a)g_2(bb^*)\\
 &&~~~~~+s^2f_2(a^*a)g_1(b^*b)
 +s^{-2}f_1(aa^*)g_2(bb^*)\big]^{1/2}.
 \ee
Now taking the infimum over all $s>0$, we finally get (\ref{G3})
with $K_3\le 2K_2$.

 iii) $\Rightarrow$ iv). This is a successive use of
the Cauchy- Schwarz and H\"older inequalities. Indeed, by
(\ref{G3}) we have
 \be
 \big|\sum_ku(a_k, b_k)\big|
 &\le& K_3\Big[\sum_k\big(f_1(a_ka_k^*)g_1(b_k^*b_k)\big)^{1/2}
 +\sum_k \big(f_2(\mu_ka_k^*a_k)g_2(\mu_k^{-1}b_kb_k^*)
 \big)^{1/2}\Big]\\
 &\le& K_3\Big[\big(\sum_kf_1(a_ka_k^*)\big)^{1/2}\,
 \big(\sum_kg_1(b_k^*b_k)\big)^{1/2}\\
 &&~~+\big(\sum_kf_2(\mu_ka_k^*a_k)\big)^{1/2}\,
 \big(\sum_kg_2(\mu_k^{-1}b_kb_k^*)\big)^{1/2}\Big]\\
 &\le& K_3\Big[\Big\|(\sum_k a_k a_k^*)^{1/2}\Big\|_p\;
 \Big\|(\sum_kb_k^*b_k)^{1/2}\Big\|_p\\
 &&~~+\Big\|(\sum_k\mu_ka_k^*a_k)^{1/2}\Big\|_p\;
 \Big\|(\sum_k\mu_k^{-1}b_kb_k^*)^{1/2}\Big\|_p\Big].
 \ee
This is (\ref{G4}) with $K_4\le K_3$. It is easy to see that iv)
$\Rightarrow$ ii) with $K_2\le K_4$ (although we will not need
this). Therefore, ii)$\Leftrightarrow$ iii) $\Leftrightarrow$ iv).

\medskip

iv) $\Rightarrow$ v). Put
 $$X=\big(E\otimes_{h_p}F\big)\oplus_1\big(F\otimes_{h_p}E\big)
 \quad \mbox{and}
 \quad Y=\{(x,\ ^tx)\ :\ x\in E\otimes F\}.$$
Then by Proposition \ref{hpq}, ii),  (\ref{G4}) implies that $u$
defines a continuous linear functional on $Y$ with norm $\le K_4$.
Hence $u$ extends to a continuous functional $\f$ on $X$ with the
same norm. Then $\f$ can be decomposed as $\f=\f_1+\f_2$ such that
 $$\f_1\in (E\otimes_{h_p}F)^*,\ \f_2\in
 (F\otimes_{h_p}E)^* \quad\mbox{and}\quad
 \max\{\|\f_1\|, \|\f_2\|\}=\|\f\|\le K_4.$$
Going back to bilinear forms on $E\times F$ and using Theorem
\ref{dual of hp}, we see that $\f_1$ and $\f_2$ define
respectively bilinear forms $u_1$ and $u_2$ such that
 $$\widetilde u_1\in\G_{R_p}(E, F^*),\
 \widetilde {^tu_2}\in\G_{R_p}(F^*, E)\quad\mbox{and}\quad
 \g_{R_p}(\widetilde u_1)=\|\f_1\|,\
 \g_{R_p}(\widetilde {^tu_2})=\|\f_2\|.$$
Then clearly,  $u=u_1+u_2$ yields the desired decomposition.

\medskip

 v) $\Rightarrow$ vi). Suppose $u_1$ and $u_2$ are as in v). Then
there are a Hilbert space $H$ and  c.b. factorizations
$\displaystyle E\mathop{\longrightarrow}^{\a_1}
H_p^r\mathop{\longrightarrow}^{\b_1} F^*$ for $\wt u_1$ and
$\displaystyle E\mathop{\longrightarrow}^{\a_2}
H_p^c\mathop{\longrightarrow}^{\b_2} F^*$ for $\wt u_2$. Define
$\a: E\to H_p^r\oplus_p H_p^c$ by $\a(a)=(\a_1(a), \a_2(a))$ and
$\b:H_p^r\oplus_p H_p^c\to F^*$ by $\b(\xi,
\eta)=\b_1(\xi)+\b_2(\eta)$. Then $\wt u=\b\a$ and
 $$\|\a\|_{cb}\le (\|\a_1\|_{cb}+\|\a_2\|_{cb})^{1/2}\ ,
 \quad
 \|\b\|_{cb}\le (\|\b_1\|_{cb}+\|\b_2\|_{cb})^{1/2}\ .$$
It then follows that $\widetilde u\in \G_{R_p\oplus_p C_p}$ and
 $$\g_{R_p\oplus_p C_p}(\widetilde u)\le
 \max\big\{\g_{R_p}(\widetilde u_1), \g_{C_p}(\widetilde u_2)
  \big\}.$$

 vi) $\Rightarrow$ i). This is evident.

 Therefore, we have proved that all assertions in Theorem \ref{G} are
equivalent. The last part of the theorem is clear from the
preceding proof. \cqd

\medskip

\n{\bf Remark.} The previous proof also works for the case of
$p=\8$ with the additional assumption that both  $E$ and $F$ are
exact. The place where we need this assumption is only the
implication i)$\Rightarrow$ii), for which we have to use Corollary
\ref{multiplication map}, ii).

\medskip

We end this section with an alternate direct proof of the
implication iii)$\Rightarrow$v) in Theorem \ref{G}. This is a
special case of the following result due to Pisier, which has
independent interest.

\begin{prop}\label{decomposition}
 Let $E, F$ be vector spaces and  $H_i, K_i$ Hilbert spaces
$(i=1, 2)$. Let $I_i: E\to H_i$ and $J_i: F\to K_i$ be linear
maps. Assume a bilinear form $u: E\times F\to \comp$ satisfies
 \beq\label{decomposition1}
 |u(a, b)|\le \|I_1(a)\|_{H_1}\|J_1(b)\|_{K_1}
  +\|I_2(a)\|_{H_2}\|J_2(b)\|_{K_2}\ ,
  \quad (a, b)\in E\times F.
 \eeq
Then $u$ can be decomposed as a sum of two bilinear forms
$u=u_1+u_2$ such that
 $$|u_1(a, b)|\le \|I_1(a)\|_{H_1}\|J_1(b)\|_{K_1}
 \quad\mbox{and}\quad
 |u_2(a, b)|\le\|I_2(a)\|_{H_2}\|J_2(b)\|_{K_2}\ ,
 \quad (a, b)\in E\times F.$$
\end{prop}

\pf On the vector space $E\ot F$ we introduce the following
semi-norm. For $x\in E\ot F$ define
 $$\|x\|_1=
 \inf\Big\{\big(\sum_k \|I_1(a_k)\|_{H_1}^2\big)^{1/2}
 \big(\sum_k \|J_1(b_k)\|_{K_1}^2\big)^{1/2}\Big\},$$
where the infimum runs over all decompositions of $x$ as $x=\sum_k
a_k\otimes b_k$. Similarly, we define a semi-norm $\|\ \|_2$ by
using $H_2$ and $K_2$. Setting
 $$a=(a_1, ... , a_n)\quad\mbox{and}\quad b=\ ^t(b_1, ... , b_n),$$
 we can rewrite $x=\sum_k a_k\otimes b_k$ as $x=a \odot b$. Then
 $$\|x\|_1=
 \inf\left\{
 \|I_{\ell_2}\otimes I_1(a)\|_{\ell_2(H_1)}
 \|I_{\ell_2}\otimes J_1(b)\|_{\ell_2(K_1)}\ ,\ x=a\odot b\right\}.$$
Since the semi-norm $\|I_{\ell_2}\otimes I_1(a)\|_{\ell_2(H_1)}$
is given by a quadratic form, for any operator  $\a\in B(\ell_2)$
 $$\|I_{\ell_2}\ot I_1(a\a)\|_{\ell_2(H_1)}\le
 \|I_{\ell_2}\ot I_1(a)\|_{\ell_2(H_1)}\ \|\a\|,$$
 $$\|I_{\ell_2}\otimes J_1(\a b)\|_{\ell_2(K_1)}\le
 \|\a\|\ \|I_{\ell_2}\otimes J_1(b)\|_{\ell_2(K_1)} .$$
From this observation and using the same argument as in the proof
of \cite[Proposition~1.7] {pisshlyak},  we can deduce the
following lemma.

\begin{lem}\label{decomposition2}
 For any $x\in E\ot F$ there is a decomposition
$x=\sum_{k=1}^n a_k\otimes b_k$ and positive numbers $\l_1, ... ,
\l_k$ such that $a_1, ..., a_n$ (resp. $b_1, ... , b_n$) are
linearly independent and such that
 \be
 \|x\|_1&=&\big(\sum_k \|I_1(a_k)\|_{H_1}^2\big)^{1/2}
 \big(\sum_k \|J_1(b_k)\|_{K_1}^2\big)^{1/2},\\
 \|x\|_2&=&\big(\sum_k \l_k\|I_2(a_k)\|_{H_2}^2\big)^{1/2}
 \big(\sum_k \l_k^{-1}\|J_2(b_k)\|_{K_2}^2\big)^{1/2}.
 \ee
\end{lem}

This lemma allows us to finish the proof of the proposition.
Indeed, let
 $$X=(E\ot F,\ \|\ \|_1) \op_1 (E\ot F,\ \|\ \|_2)\quad\mbox{and}
 \quad Y=\{(x, x)\ :\ x\in E\ot F\}\subset X.$$
By (\ref{decomposition1}) and Lemma \ref{decomposition2},
considered as a linear functional on $E\ot F$ (and so on $Y$ too),
$u$ is continuous and of norm $\le 1$ with respect to the
semi-norm of $Y$. Therefore, by the Hahn-Banach theorem, $u$
extends to a contractive functional $\wh u$ on $X$. Write
 $$\wh u(x, y)=\wh u(x, 0) +\wh u(0, y)\mathop{=}^{\rm def}
 u_1(x)+u_2(y),\quad x, y\in E\ot F.$$
Considered back to bilinear forms on $E\times F$, $u_1$ and $u_2$
give the required decomposition of $u$. \cqd

\medskip

\n{\bf Remark.} As the reader can see, the above proof is similar
to that of the implication iv) $\Rightarrow$ v) of Theorem
\ref{G}. However, it has an advantage that  the original
functionals $f_i$ and $g_i$ in (\ref{G3}) can be used for $u_1$
and $u_2$ respectively in Theorem \ref{G}, v); see Theorem
\ref{dual of hp} for the existence of such functionals for $u_1$
and $u_2$. (Concerning this point see a remark in
\cite[p.189]{pisshlyak}.)

\section{Proof  of Theorem \ref{LG}}\label{pfLG}

Now we pass to the proof of Theorem \ref{LG}. We will need the
following result, which is a generalization of \cite[Theorem 8.4]
{pis-OH}.

\begin{prop}\label{SpCq}
Let $2<p<\8, 0\le\th\le1$ and $\frac{1}{q}=\frac{1-\th}{p}+
\frac{\th}{p'}$. Let $(e_k)$ denote the canonical basis of $C_q$.
\begin{enumerate}[{\rm i)}]
\item For any finite sequence $(x_k)\subset S_p$ we have
 \beq\label{SpCq1}
 \|\sum_kx_k\ot e_k\|_{S_p[C_q]}
 =\sup\Big\{\big(\sum_k\|\a x_k\b\|_2^2\big)^{1/2}\Big\},
 \eeq
where the supremum runs over all $\a$ and $\b$ respectively in the
unit balls of $S_{2r\th^{-1}}$ and $S_{2r(1-\th)^{-1}}$, $r$ being
the conjugate index of $p/2$. Moreover, the supremum can be
restricted to all $\a$ and $\b$ in the positive parts of these
unit balls.
\item Let $H$ be a Hilbert space and $\U$ a free ultrafilter. Let
$B(H)_{\U}$ be the associated ultrapower von Neumann algebra. Then
for any finite sequence $(x_k)\subset L_p(B(H)_{\U})$
 \beq\label{SpCq2}
 \|\sum_kx_k\ot e_k\|_{L_p[B(H)_{\U};C_q]}
 =\sup\Big\{\big(\sum_k\|\a x_k\b\|_2^2\big)^{1/2}\Big\},
 \eeq
where the supremum runs over all $\a$ and $\b$ respectively in the
unit balls of $L_{2r\th^{-1}}(B(H)_{\U})$ and
$L_{2r(1-\th)^{-1}}(B(H)_{\U})$; again the supremum can be
restricted to all $\a$ and $\b$ in the positive parts of these
unit balls.
\end{enumerate}
\end{prop}

\pf i) Let $x=\sum_kx_k\ot e_k\in S_p[C_q]$.  By (\ref{Sq via
Sp}),
 $$\|x\|_{S_p[C_q]}=\sup\big\{\|axb\|_{S_q[C_q]}\;:\; a, b\in
 S_{2s},\ \|a\|_{2s}\le 1,  \|b\|_{2s}\le 1\big\},$$
where $s$ is determined by $\frac{1}q=\frac{1}p+\frac{1}s$.
However, $S_q[C_q]=C_q[S_q]$. Thus by (\ref{col Lp}),
 $$\|axb\|_{S_q[C_q]}^2=\big\|\sum_k
 b^*x_k^*a^*ax_kb\big\|_{\frac{q}2}\,.$$
Assume $q\ge2$. It then follows that
 $$\|x\|_{S_p[C_q]}^2=\sup\big\{\sum_k
 \Tr(c^*b^*x_k^*a^*ax_kbc): a, b\in
 S_{2s},c\in S_{2t}, \|a\|_{2s}\le 1,
 \|b\|_{2s}\le 1,\|c\|_{2t}\le 1 \big\},$$
where $t$ is the index conjugate to $q/2$. Set $\a=a$ and $\b=bc$.
Note that $\a\in S_{2r\th^{-1}}$ and $\b\in S_{2r(1-\th)^{-1}}$.
We then deduce (\ref{SpCq1}) in the case of $q\ge2$. The case of
$q\le2$ can be done similarly by using the identification
$C_q\cong R_{q'}$ (see (\ref{dual Cp-Rp})).

ii) It suffices to prove (\ref{SpCq2}) for any $C_q^n$ instead of
$C_q$. Then by definition, $L_p[B(H)_{\U}; C_q^n]$ is the
ultrapower  $\big(S_p[H; C_q^n]\big)^{\U}$. On the other hand,
$L_p(B(H)_{\U})$ is also the ultrapower  $\big(S_p(H)\big)^{\U}$
(Raynaud's theorem). Recall that again by [Ra], ultraproduct
preserves all algebraic structures on noncommutative $L_p$-spaces,
in particular, the product. With the help of all these, we can
easily deduce (\ref{SpCq2}) from (\ref{SpCq1}).\cqd

\medskip

\n{\bf Remarks.} i) We are grateful to the referee for the short
proof of part i) above, which is much simpler than our original
one.

ii) One can show that Proposition \ref{SpCq}, ii) holds for any
QWEP von Neumann algebra $M$ in place of $B(H)_{\U}$.

iii) Proposition \ref{SpCq} yields a simple description of the
norm in the complex interpolation space $(C_p[L_p(M)],\;
R_p[L_p(M)])_\th$ when $M$ is $B(H)$ or an ultrapower of $B(H)$
($p\ge 2$). More generally, one can describe the norm of
$(C_p[L_p(M)],\; R_p[L_p(M)])_\th$ for any von Neumann algebra $M$
by a formula like (\ref{SpCq1}) in the case $p\ge 2$, and by a
similar dual formula in the case $p<2$. This will be pursued
elsewhere.

\medskip

\n{\it Proof of Theorem 0.2.} Without loss of generality, we can
assume $H$ separable and infinite dimensional. Thus $H_q^c=C_q$

i) $\Rightarrow$ ii). Assume $\|u\|_{cb}\le 1$. As for the proof
of Theorem \ref{G}, the noncommutative Khintchine inequality in
section~3 will be the key ingredient for the present proof too.
Fix a positive sequence $(\l_k)$. We maintain the notations
introduced at the beginning of the proof of Theorem \ref{G}, but
with $g_k$ being now defined by
 $$g_k=\l_k^{\th}\;\el(e_k)+ \l_k^{-(1-\th)}\,\el^*(e_{-k}).$$
The $g_{k, p}$ are defined by (\ref{circular p}) with the $g_k$
above. Then by Proposition \ref{cb extension},
 $$\|I_{L_p(\G)}\ot u\|\le 1.$$
Thus by Theorem \ref{khintchine}, for any finite sequence
$(a_k)\subset E$
 \be
 &&\Big\|\sum_k u(a_k)\ot g_{k, p}\Big\|_{L_p[\G; C_q]}
 \le\Big\|\sum_k a_k\ot g_{k, p}\Big\|_{L_p[\G; E]}\\
 &&~~\le B_p \max\Big\{
 \Big\|\big(\sum_k\l_k^{2\th(1-\frac{2}{p})}a_k^*a_k
 \big)^{\frac{1}{2}}\Big\|_{p}\ ,
 \Big\|\big(\sum_k\l_k^{-2(1-\th)(1-\frac{2}{p})}a_ka_k^*
 \big)^{\frac{1}{2}}\Big\|_{p}\Big\}.
 \ee
Therefore, to prove (\ref{LG2}) we must show that for any finite
sequence $(z_k)\subset C_q$
 \beq\label{i to ii}
 \big(\sum_k\|z_k\|^2\big)^{\frac{1}{2}}\le
 \big\|\sum_k z_k\ot g_{k, p}\big\|_{L_p[\G; C_q]}\;.
 \eeq
To this end we can clearly assume all $z_k$ are finitely
supported, and so $C_q$ can be replaced by a $C_q^n$. Write $z_k$
in the canonical basis of $C_q^n$:
 $$z_k=\sum_{j=1}^n z_{k,j}\,e_j\;.$$
Then
 $$\sum_k z_k\ot g_{k,p}=\sum_{j=1}^n e_j\ot x_j,$$
where
 $$x_j=\sum_k z_{k,j}\,g_{k, p}\in L_p(\G).$$
Now since $\G$ is QWEP, $\G$ is the image of a normal conditional
expectation $\Phi$ on some ultrapower von Neumann algebra
$B(H)_{\U}$. Then by definition, $L_p[\G; C_q^n]$ is a subspace of
$L_p[B(H)_{\U}; C_q^n]$. By Proposition \ref{SpCq}, ii), for any
$\a$ and $\b$  in the unit balls of $L_{2r\th^{-1}}(B(H)_{\U})$
and $L_{2r(1-\th)^{-1}}(B(H)_{\U})$, respectively, we have
 \beq \label{i to iibis}
 \big(\sum_j\|\a x_j\b\|_2^2\big)^{1/2}\le
 \|\sum_je_j\ot x_j\|_{L_p[\G;C^n_q]}\;.
 \eeq
In particular, this is true for $\a=D^{\th/2r}$ and
$\b=D^{(1-\th)/2r}$. For this choice of $\a$ and $\b$, we have
 \be
 \|\a x_j\b\|_2^2
 &=&\|\a x_j\b\|_{L_2(\G)}^2
 =\tr(\b^*x_j^*\a^* \a x_j\b)\\
 &=&\sum_{k, k'}\overline{z_{k,j}}\,z_{k',j}\,\tr\big[g_k^*D^\th
 g_{k'}D^{1-\th}\big]\\
 &=&\sum_k|z_{k,j}|^2\quad\mbox{by}\; (\ref{circular orthogonal}).
 \ee
Summing up over all $j$ and using (\ref{i to iibis}), we get
(\ref{i to ii}). Therefore, (\ref{LG2}) is proved.

\medskip

ii) $\Rightarrow$ iii). This is a standard application of the
Hahn-Banach theorem as in the proof of ii) $\Rightarrow$ iii) in
Theorem \ref{G}. Conversely, it is trivial that iii) $\Rightarrow$
 ii).

iii) $\Rightarrow$ i).  Let $a=(a_{i\,j})\in S^n_p[E]$ be a unit
element. We have to show
 $$\|I_{S^n_p}\ot u(a)\|_{S^n_p[C_q]}\le K.$$
To this end we use again Proposition \ref{SpCq}. Let $u_k$ be the
k-th component of $u$ relative to the canonical basis of $C_q$.
Set $x_k=\big(u_k(a_{ij})\big)_{1\le i, j\le n}\in S_p^n$. Then
 $$I_{S^n_p}\ot u(a)=\sum_k x_k\ot e_k.$$
Let $\a$ (resp. $\b$) be a positive matrix in the unit ball of
$S_{2r\th^{-1}}$ (resp. $S_{2r(1-\th)^{-1}}$). We are going  to
estimate
 $\sum_k\|\a x_k\b\|_2^2.$
In virtue of the invariance of the norm of $S^n_p[E]$ by
multiplication from left and right by unitary matrices, and
changing the matrix $a$ if necessary, we can assume both $\a$ and
$\b$ are diagonal. Let $\a_i$ and $\b_i$ be respectively their
diagonal entries. Then by (\ref{LG3}) and the H\"older
inequalities
 \be
 \sum_k\|\a x_k\b\|_2^2
 &=&\sum_k\sum_{i, j}\a_i^2\b_j^2\,|x_k(i, j)|^2
 =\sum_{i\, j}\a_i^2\b_j^2\,\|u(a_{ij})\|^2\\
 &\le& K^2\, \sum_{i,
 j}\a_i^2\b_j^2\,[f(a_{ij}^*a_{ij})]^{1-\th}\,
 [g(a_{ij}a_{ij}^*)]^{\th}\\
 &\le& K^2\, \Big(\sum_{i,j}\b_j^{2(1-\th)^{-1}}\;f(a_{ij}^*a_{ij})
 \Big)^{1-\th}\;\Big(\sum_{i, j}\a_i^{2\th^{-1}} g(a_{ij}a_{ij}^*)
 \Big)^{\th}\\
 &\le& K^2\, \Big(\sum_{j}\b_j^{2(1-\th)^{-1}}\Big\|\sum_{i}a_{i j}^*
 a_{ij}\Big\|_{L_{\frac{p}{2}}(M)}\Big)^{1-\th}
 \Big(\sum_{i}\a_i^{2\th^{-1}}\Big\|
 \sum_{j}a_{ij}a_{ij}^*\Big\|_{L_{\frac{p}{2}}(M)}\Big)^{\th}\\
 &\le& K^2\,\big\|\b\big\|^{\frac{1-\th}{r}}_{2r(1-\th)^{-1}}
 \Big(\sum_{j}\Big\|\sum_{i}a_{ij}^*a_{ij}
 \Big\|^{\frac{p}{2}}_{L_{\frac{p}{2}}(M)}
 \Big)^{\frac{2(1-\th)}{p}}\\
 &&~~\bullet\big\|\a\big\|^{\frac{\th}{r}}_{2r\th^{-1}}
 \Big(\sum_{i}
 \Big\|\sum_{j}a_{ij}a_{i j}^*\Big\|^{\frac{p}{2}}_{L_{\frac{p}{2}}(M)}
 \Big)^{\frac{2\th}{p}}\\
 &\le& K^2\, \|a\|^2_{S_p[E]}\; ,
 \ee
where for the last inequality we have used the following
elementary fact that for any von Neumann algebra $M$ and any $a\in
S_p[L_p(M)]$ with $2\le p\le\8$
 \be
 && \big(\sum_{j}\big\|(\sum_{i}a_{ij}^*a_{i
 j})^{1/2}\big\|^{p}_{L_{p}(M)}\big)^{1/ p}\le \|a\|_{S_p[L_p(M)]},\\
 &&\big(\sum_{i}\big\|(\sum_{j}a_{i\,j}a_{i\, j}^*)^{1/2}
 \big\|^{p}_{L_{p}(M)}\big)^{1/ p}\le \|a\|_{S_p[L_p(M)]}.
 \ee
Therefore, $\|I_{S^n_p}\ot u(a)\|_{S^n_p[C_q]}\le K$ for any $n\ge
1$. Then by Lemma \ref{cb via Sp}, $\|u\|_{cb}\le K$. Therefore,
we have proved Theorem \ref{LG}. \cqd

\section{Applications}\label{appli}

In this section we present some applications of Theorems 0.1 and
0.2. We first give a factorization for maps satisfying one of the
conditions in Theorem \ref{LG} through a real interpolation space
of parameters $(\th, 1)$. To this end let $f$ and $g$ be two
positive unit functionals  on $L_{p/2}(M)$. Let $K_{f}$ be the
Hilbert space obtained from $E$ relative to the semi-scalar
product $\la a, b\ra=f(a^*b)$ (see Remark \ref{iE} and  the proof
of iv) $\Rightarrow$ i) of Proposition \ref{to Rp}). Similarly,
the semi-scalar product $(a, b)\mapsto g(ba^*)$ yields another
Hilbert space $K_{g}$. Let $i_f: E\to K_f$ (resp. $i_g: E\to K_g$)
be the natural inclusion. Note that both $i_f$ and $i_g$ are
injective and of dense range. It follows that $i_f^*: K_f^*\to
E^*$ and $i_g^*: K_g^*\to E^*$ are injective. This allows us to
regard $(K_f^*\,,\; K_g^*\,)$, and so $(K_f\,,\; K_g)$ as
compatible couples of Hilbert spaces. Under this compatibility,
$i_f$ and $i_g$ are the same map, denoted by $i_{f, g}$ below.

Now we equip $K_f$ (resp. $K_g$) with the operator space structure
of $K_{f,p}^c$ (resp. $K_{g,p}^r$), and consider real
interpolation space $(K_{f,p}^c\,,\; K_{g,p}^r\,)_{\th, 1}$ (see
\cite{xu-int} for the real interpolation theory in the category of
operator spaces). By Remark \ref{iE}, $i_f: E\to K_{f,p}^c$ and
$i_g: E\to K_{g,p}^r$ are completely contractive, so by
interpolation they induce a completely contractive map $i_{f, g}:
E\to (K_{f,p}^c\,,\; K_{g,p}^r\,)_{\th, 1}$.

\begin{cor}\label{LGfac} Let $E, p, q, \th$ be as in Theorem
\ref{LG}. Then a map $u: E\to H_q^c$  is c.b. iff there are two
positive unit functionals $f$ and $g$ on $L_{p/2}(M)$ such that
$u$ admits a factorization of the following form
 \beq\label{LGfac1}
 \xymatrix{E \ar[rr]^u \ar[rd]_{i_{f, g}}& &  H_q^c
 \\ & (K_{f,p}^c\,,\; K_{g,p}^r)_{\th, 1}\ar[ur]_{\wh u}\;, &}
 \eeq
where $\wh u$ is a \underline{bounded} map and $\|\wh u\|\le K$.
Moreover, the smallest of  such constants $K$ is universally
equivalent to $\|u\|_{cb}$.
\end{cor}

\pf Let $u: E\to H_q^c$ be c.b.. Then by Theorem \ref{LG}, $u$
satisfies (\ref{LG3}). By the previous discussion, we clearly have
the required factorization. Conversely, if $u$ admits such a
factorization, then we have (\ref{LG3}), and so $u$ is c.b.. \cqd

\medskip

\n{\bf Remarks.} i) Using \cite{pis-Rp}, we can get a
factorization similar to that in Corollary \ref{LGfac} in the case
of $p=\8$ with the additional assumption that either $E\subset M$
is exact or $E=M$.

ii) We do not know whether $\wh u$ in Corollary \ref{LGfac} can be
chosen to be c.b..

\medskip

Let us isolate out the special case of $\th=1/2$ in Theorem
\ref{LG} and Corollary \ref{LGfac} because of the particular
importance of $OH$.

\begin{cor}\label{LGOH}
 Let $E\subset L_p(M)$ be a subspace with $2<p\le\8$. In the
case of $p=\8$ we assume in addition that $E$ is either exact with
${\rm ex}(E)=1$ or $E=M$. Then for any map $u: E\to OH(I)$ $($with
$I$ an index set$)$ the following assertions are equivalent:
\begin{enumerate}[{\rm i)}]
\item  $u$ is c.b..
\item There is a constant $K$ such that for all finite
sequences $(a_k)\subset E$ and  $(\mu_k)\subset\real_+$
 \beq\label{LGOH2}
 \sum_k\|u(a_k)\|^2\le {\frac{K^2}{2}} \Big[
 \big\|\sum_k\mu_ka_k^*a_k\big\|_{p/2}+
 \big\|\sum_k\mu_k^{-1}a_ka_k^*\big\|_{p/2}\Big].
 \eeq
\item There are two positive unit functionals $f, g$ on
$L_{p/2}(M)$ such that
 \beq\label{LGOH3}
 \|u(a)\|\le K\, \big(f(a^*a)\big)^{1/4}
 \big(g(aa^*)\big)^{1/4}\ ,\quad a\in E.
 \eeq
\item There are two positive unit functionals $f, g$ on
$L_{p/2}(M)$ such that $u$ admits a factorization of the form
(\ref{LGfac1})
 with $\|\wh u\|\le K'$.
\end{enumerate}
Moreover, the best constants $K$ and $K'$ are universally
equivalent to $\|u\|_{cb}$.
\end{cor}

This is the little Grothendieck theorem for noncommutative
$L_p$-spaces in the category of operator spaces. The case of
$p=\8$ goes back to \cite{pisshlyak}. Compare (\ref{LGOH3}) with
(\ref{piquard}): the arithmetic mean in (\ref{piquard}) is
replaced by the geometric mean in (\ref{LGOH3}).

\begin{cor}\label{extension}
  In the situation of Theorem \ref{G}, if
$u: E\times F\to\comp$ is j.c.b., then $u$ admits an extension $U:
L_p(M)\times L_p(M)\to\comp$ with $\|U\|_{jcb}\le c\|u\|_{jcb}$,
where $c$ is a universal constant.
\end{cor}

\pf  We use the decomposition $u=u_1+u_2$ in Theorem \ref{G}, v).
Regarding $u_1, u_2$ as linear functionals on $E\ot F$, we see
that $u_1$ and $ ^tu_2$ satisfy iv) of Theorem \ref{dual of hp}.
Therefore, $u_1$ (resp. $u_2$) admits an extension  $U_1$ (resp.
$U_2$) on $L_p(M)\ot_{h_p}L_p(M)$ (resp. $L_p(M)\ot_{h_p}L_p(M)$).
Then $U= U_1+U_2$ is the required extension of $u$.\cqd

\medskip

We say that a c.b. map $T: E\to F$ between two operator spaces has
the {\it completely bounded approximation property} (CBAP in
short) if there are a constant $\l$ and  a net $(T_i)$ of finite
rank maps from $E$ to $F$ such that $T_i$ converges to $T$ in the
point-norm topology and $\sup_i\|T_i\|_{cb}\le\l \|T\|_{cb}$. In
this case, we also say that $T$ has the $\l$-CBAP if we want to
emphasize the constant $\l$. Note that $E$ has the CBAP iff the
identity of $E$ does.  We recall the open problem in
\cite{pisshlyak} whether any c.b. map from a C*-algebra to the
dual of a C*-algebra has automatically the CBAP. However, the
corresponding problem in the $L_p$-space case is easily solved by
virtue of Theorem \ref{G}.

\begin{cor}
 Let $E, F$ be as in Theorem \ref{G}. Then any map $T\in
CB(E, F^*)$ has the $\l$-CBAP with $\l$ a universal constant.
\end{cor}

\pf By Theorem \ref{G}, $T$ belongs to $\G_{R_p\op C_p}(E, F^*)$.
It remains to note that $R_p\op C_p$ has the 1-CBAP. \cqd

\begin{cor} {\rm i)} Let $E$ be an operator space. If both
$E$ and $E^*$ are completely isomorphic to subspaces of a
noncommutative $L_p(M)$ with $1<p<2$, then $E$ is completely
isomorphic to a quotient of a subspace of $H_p^c\op_p K_p^r$ for
some Hilbert spaces $H$ and $K$.\hfill\break\indent
 {\rm ii)} If we assume in addition that the completely isomorphic
copies of $E$ and $E^*$ are completely complemented in $L_p(M)$,
then $E$ is completely isomorphic to $H_p^c\op_p K_p^r$.
\end{cor}

\pf Based on Theorem \ref{G}, the proof of this corollary is the
same as that of Corollaries 3.1 and 3.3 in \cite{pisshlyak} (which
corresponds to the case $p=1$ with an additional assumption on
$E$), so we omit it.\cqd

\medskip

\n{\bf Remarks.} i) Let $H$  be a Hilbert space  and $1<p<2$. It
is proved in \cite{xu-emb} that a quotient of a subspace of
$H_p^c\op_p H_p^r$ is completely isomorphic to a subspace of a
noncommutative $L_p$.

ii) Moreover, for any  $p<q\le 2$, $H_q^c$ is a quotient of a
subspace of $H_p^c\op_p H_p^r$. Consequently, any quotient of a
subspace of $H_q^c\op_q H_q^r$ is also completely isomorphic to a
subspace of a noncommutative $L_p$ ($1\le p<q\le 2$). See
\cite{xu-emb}  for more details.

\medskip

We end this section with an application to Schur multipliers.  Let
$\f$ be a function on $\nat\times \nat$. We recall that $\f$ is a
Schur multiplier from $S_p$ to $S_q$ if the map $M_\f: x\mapsto
(\f(i,j)x(i,j))$ defined for finite matrices $x$ extends to a
bounded map from $S_p$ to $S_q$ (which is still denoted by
$M_\f$). (Note that we change slightly the matrix notation by
regarding a matrix as a function on $\nat\times \nat$ too.) If
$M_\f$ is c.b., we say that $\f$ is a c.b. Schur multiplier from
$S_p$ to $S_q$. More generally, if $\La\subset \nat\times \nat$ is
a subset, we denote by $S_p^\La\subset S_p$ the subspace of all
$x\in S_p$ which vanish outside $\La$. A function $\f:
\La\to\comp$ is called a (c.b.) Schur multiplier from $S_p^\La$ to
$S_q$ if $x\mapsto (\f(i,j)x(i,j))$ extends to a bounded (c.b.)
map from $S_p^\La$ to $S_q$.

\medskip

$\ell_r(\ell_\8)$ is the space of all complex functions $\f$ on
$\nat\times \nat$ such that
 $$\|\f\|_{\ell_r(\ell_\8)}=\big(\sum_i\sup_j|\f(i,j)|^r\big)^{1/r}
 <\8.$$
Set
 $ ^t\ell_r(\ell_\8)=\{\f\ :\ ^t\f\in \ell_r(\ell_\8)\},$
where $ ^t\f(i,j)=\f(j,i)$.

\medskip

The following, except ${\rm i}'$) and ${\rm ii}'$),  is again the
$L_p$-space analogue of the corresponding results in
\cite{pisshlyak}, which correspond  to the case where $q=1$ and
$p=\8$.

\begin{cor}\label{schur}
 {\rm i}) Let $1\le q\le 2\le p\le\8$ and
$r=\frac{pq}{p-q}$. Then $\f$ is a Schur multiplier from $S_p$ to
$S_q$ iff $\f\in \ell_r(\ell_\8)+\, ^t\ell_r(\ell_\8)$, i.e. iff
$\f$ admits a decomposition
 $$\f=\f_1+\f_2\quad\mbox{with}\quad \f_1, \;
 ^t\f_2\in \ell_r(\ell_\8).$$
Moreover,
 $$\|M_\f\|\approx
 \inf\big\{\|\f_1\|_{\ell_r(\ell_\8)} + \|\, ^t\f_2\|_{\ell_r(\ell_\8)}\ :\
 \f=\f_1+\f_2,\; \f_1,\; ^t\f_2\in \ell_r(\ell_\8)\big\}.$$
 \indent ${\rm i}'$) Assume in addition $p<\8$ and
 $\La\subset \nat\times \nat$. Then every Schur multiplier from $S^\La_p$ to
$S_q$ extends to a Schur multiplier from $S_p$ to
$S_q$.\hfill\break \indent
 {\rm ii}) Let $2\le p\le\8$ and $r$ be the
conjugate index of $p/2$. Then $\f$ is a c.b. Schur multiplier
from $S_p$ to $S_{p'}$ iff there are $\a,\b\in \ell_{2r}$ such
that $|\f_{ij}|\le |\a_i||\b_j|$ for all $i,j$, or equivalently,
iff $\f$ admits a factorization
 $$|\f|=(|\f_1|\, |\f_2|)^{1/2}
 \quad\mbox{with}\quad \f_1, \; ^t\f_2\in \ell_r(\ell_\8).$$
Moreover, in this case
 $$\|M_\f\|_{cb}\approx \inf\big\{\|\a\|_{2r}\|\b\|_{2r}\ :\
 |\f_{ij}|\le |\a_i||\b_j|, \ \a,\b\in \ell_{2r}\big\}.$$
 \indent ${\rm ii}'$) With the same $p$ and $r$ as in {\rm ii)},
let $\La\subset \nat\times \nat$. Then every c.b. Schur multiplier
from $S^\La_p$ to $S_{p'}$ extends to a c.b. Schur multiplier from
$S_p$ to $S_{p'}$.
\end{cor}

\pf The proof is similar to those of Theorems 4.1 and 4.2 in
\cite{pisshlyak}. Thus we will be very brief.

i) The case $p=\8$ and $q=1$ corresponds to \cite[Theorem 4.1]
{pisshlyak}. Let $\f$ be a Schur multiplier from $S_p$ to $S_q$
with $\|M_\f\|\le 1$. Let $u: S_p\times S_{q'}\to \comp$ be the
bilinear form defined by $M_\f$:
 $$u(x, y)=\sum_{i,j}\f(i,j)x(i,j)y(i,j),\quad x\in S_p,\; y\in
 S_{q'}.$$
Then by (\ref{piquard1}), $u$ can be decomposed as
$u=u_1+u_2+u_3+u_4$ with
 $$|u_1(x, y)|\le K\big(f_1(xx^*)g_1(y^*y)\big)^{1/2},\quad
 |u_2(x, y)|\le K\big(f_2(x^*x)g_2(yy^*)\big)^{1/2},$$
 $$|u_3(x, y)|\le K\big(f_1(xx^*)g_2(yy^*)\big)^{1/2},\quad
 |u_4(x, y)|\le K\big(f_2(x^*x)g_1(y^*y)\big)^{1/2},$$
where $f_1, f_2$ (resp. $g_1, g_2$) are positive unit functionals
on $S_{p/2}$ (resp.  $S_{q'/2}$). By an elementary average
argument as in \cite{pisshlyak}, we can assume that  each $u_k$ is
given by a Schur multiplier as $u$, say $\f_k$, and the
functionals $f_i, g_i$ are diagonal matrices. Then we have
$\f=\f_1+\f_2+\f_3+\f_4$ with
 $$|\f_1(i,j)|\le K\,\big(f_1(i,i)g_1(j,j)\big)^{1/2},\quad
 |\f_2(i,j)|\le K\,\big(f_2(j,j)g_2(i,i)\big)^{1/2},$$
 $$|\f_3(i,j)|\le K\,\big(f_1(i,i)g_2(i,i)\big)^{1/2},\quad
 |\f_4(i,j)|\le K\,\big(f_2(j,j)g_1(j,j)\big)^{1/2}.$$
Thus $\f_3$ and $ ^t\f_4$ are in $\ell_r(\ell_\8)$. On the other
hand, $\f_1$ and $\f_2$ can be decomposed into  sums of two such
elements, i.e.  $\f_1, \f_2\in \ell_r(\ell_\8) +\,
^t\ell_r(\ell_\8)$. Hence $\f\in \ell_r(\ell_\8) +\,
^t\ell_r(\ell_\8)$.

Conversely, suppose $\f\in \ell_r(\ell_\8)$. Then for any $x\in
S_p$ and $y\in S_{q'}$
 \be
 |\sum_{i,j}\f(i,j)x(i,j)y(i,j)|
 &\le& \sum_{i}\sup_j|\f(i,j)|\,\sum_j|x(i,j)y(i,j)|\\
 &\le& \sum_{i}\sup_j|\f(i,j)|\,\big(\sum_j|x(i,j)|^2\big)^{1/2}
  \,\big(\sum_j|y(i,j)|^2\big)^{1/2}\\
 &\le& \|\f\|_{\ell_r(\ell_\8)}\|x\|_{S_p}\|y\|_{S_{q'}},
 \ee
where we have used the following elementary inequality
 $$\Big(\sum_{i}\big(\sum_j|x(i,j)|^2\big)^{p/2}\Big)^{1/p}
 \le\|x\|_{S_p}.$$
Therefore $\f$ is a Schur multiplier from $S_p$ to $S_q$.
Similarly, every matrix in $ ^t\ell_r(\ell_\8)$ is also a Schur
multiplier from $S_p$ to $S_q$. This proves part i).

${\rm i}'$) Let $\f$ be a Schur multiplier from $S^\La_p$ to
$S_q$. Since $S_p$ and $S_q$ have respectively type 2 and cotype
2, by Kwapien's theorem (cf. \cite[Corollary 3.6] {pis-cbms}),
$M_\f: S^\La_p\to S_q$ factors through a Hilbert space. Then
Maurey's extension theorem \cite{mau-pro} implies that $M_\f$
admits an extension $T: S_p\to S_q$. Averaging $T$ over the group
of all unitary diagonal matrices in $B(\ell_2)$, we deduce an
extension of $M_\f$ which is again a Schur multiplier.

ii) This part is proved in a way similar to that of i); the only
difference is that this time instead of (\ref{piquard1}) we use
Theorem \ref{G}. We omit the details.

${\rm ii}'$) Assume $M_\f: S^\La_p\to S_{p'}$ is c.b.. Then $M_\f$
has a c.b. extension from $S_p\to S_{p'}$. This follows from
Corollary 0.6 of \cite{pisshlyak} in the case $p=\8$, and from
Corollary \ref{extension}  for $p<\8$. Then as above for ${\rm
i}'$), we get a Schur extension of $\f$. \cqd

\medskip

The second part of Corollary \ref{schur} gives a characterization
of c.b. Schur multipliers from $S_p$ to $S_{p'}$, i.e. from $S_p$
to its dual. We do not know how to characterize the c.b. Schur
multipliers from $S_p$ into $S_q$ for any $1\le q\le 2\le p\le\8$,
as in the first part at the Banach space level.

\begin{pb}
 Let $1\le q\le 2\le p\le\8$ with $p\not= q'$. Characterize c.b. Schur
multipliers from $S_p$ into $S_q$, in a way similar to that in
Corollary \ref{schur}.
\end{pb}

This problem might be related to the following

\begin{pb}
 Let $M$ be a von Neumann algebra and $p, q\ge 2$ with $p\not= q$.
Let $E\subset L_p(M)$ and $F\subset L_q(M)$. Find a Grothendieck
type inequality for j.c.b. forms $u: E\times F\to\comp$.
\end{pb}

\bigskip

\n{\bf Acknowledgements.} We are grateful to Gilles Pisier for
pointing out a gap in our first proof of Theorem \ref{khintchine}
and   for allowing us to include Proposition \ref{decomposition}.
We would like to thank also Marius Junge for putting
\cite{ju-qwep} at our disposition.


\end{document}